\documentclass[11pt,a4paper, abstracton,fleqn]{scrartcl}

 \pdfoutput=1

\usepackage[latin1]{inputenc}

\usepackage{amsthm,dcolumn,graphicx,multicol,fancyhdr}
\usepackage{latexsym,subfigure}
\usepackage{sublabel}
\usepackage{rotating,amsfonts, epic,eepic,color,amssymb, amsmath, amscd, array, enumerate,afterpage}

\usepackage{hyperref}

\newtheorem{theorem}{Theorem}[section]

\newtheorem{lemma}{Lemma}[section]

\theoremstyle{definition}
\newtheorem{ex}{Example}[section]
\newtheorem{rem}{Remark}[section]

\renewenvironment{proof}[1][\negthickspace]{\textbf{Proof\thickspace#1.} }{\
\rule{0.5em}{0.5em}}

\makeatletter\@addtoreset{equation}{section}\makeatother
\makeatletter\@addtoreset{table}{section}\makeatother
\makeatletter\@addtoreset{figure}{section}\makeatother

\newcommand{\var}{\operatorname{var}}
\newcommand{\E}{\operatorname{E}}
\newcommand{\cov}{\operatorname{cov}}
\newcommand{\ls}{\leqslant}
\newcommand{\gs}{\geqslant}
\newcommand{\eps}{\epsilon}

\newcommand{\dto}{\overset{\mathcal{D}}{\longrightarrow}}

\newcommand{\wmu}{\widehat{\mu}}
\newcommand{\wm}{\widehat{m}}
\newcommand{\we}{\widetilde{e}}
\newcommand{\vth}{\vartheta}

\newcommand{\hus}{Hu\v{s}kov\'{a} }

\newcommand{\haj}{H\'{a}jek }

\begin{document}

\title{Bootstrapping confidence intervals for the change-point of time series}

\author{Marie Hu\v{s}kov\'{a}\footnote{Charles University of Prague, Department of
Statistics, Sokolovsk\'a 83, \newline
CZ\,--\,186\,75 Praha 8,
Czech Republic; \texttt{huskova@karlin.mff.cuni.cz}}, Claudia Kirch\footnote{Technical University Kaiserslautern, Departement of Mathematics, Erwin-Schr\"odinger-Stra\ss e, \newline
D--67\,653 Kaiserslautern, Germany;
\texttt{ckirch@mathematik.uni-kl.de}}}

\maketitle

\begin{abstract}
We study an AMOC time series model with an abrupt change in the mean and dependent
errors that fulfill certain mixing conditions. We obtain confidence intervals for
the unknown change-point via bootstrapping methods.\\
Precisely we use a block bootstrap of the estimated centered error sequence.
Then we reconstruct a sequence with a change in the mean using the same estimators
as before. The difference between the change-point estimator of the
resampled sequence and the one for the original sequence can be used
as an approximation of the difference between the real change-point and
 its estimator. This enables us to construct confidence intervals using
 the empirical distribution of the resampled time series.\\
A simulation study shows that the resampled confidence intervals are
 usually closer to their target levels and at the same time smaller than the asymptotic intervals.
\end{abstract}

 \vspace{5mm}
 {\bf Keywords:} confidence intervals, block bootstrap, mixing, change in mean

 \vspace{5mm}
 {\bf AMS Subject Classification 2000:} 62G09, 62G15, 60G10

 \vspace{5mm}
  {\bf Acknowledgement}: The work of the  first  author
    was partly supported by  the grants GA\v CR 201/06/0186 and
       MSM 02162839.

\section{Introduction}
Recently a number of papers has been published on possible
application of bootstrapping or permutation methods in
change-point analysis, confer Hu\v{s}kov\'{a} \cite{surveyhus} for
a recent survey. Most of these papers are concerned with obtaining
critical values for the corresponding change-point tests. Another
important issue in change-point analysis, however, is how to
obtain confidence intervals for the change-point. In this paper we
construct bootstrapping confidence intervals for the change-point in
a model with dependent data.

We consider the following At-Most-One-Change (AMOC) location model
\begin{equation}
\label{LPmodel} X(i)=\mu + d \,1_{\{i>m\}}+e(i),\quad 1\ls i \ls
n,
\end{equation}
where $m=m(n)=\lfloor n \vth\rfloor$, $d=d_n$ may depend on $n$. The errors $\{e(i),1\ls i \ls n\}$ are stationary and strong-mixing with a rate specified below, 
\begin{equation}\begin{split}\label{vor_fehler_1}
&\E e(i)=0,\quad 0<\sigma^2=\E e(i)^2<\infty,\quad    \E|e(i)|^{\nu}<\infty \text{ for some }\nu>4,\\
   &\sum_{h\gs 0}|\cov(e(0),e(h))|<\infty,
\end{split}\end{equation}
and
\begin{equation}
	\tau^2:=\var e(0)+2\sum_{h\gs 1}\cov(e(0),e(h))<\infty.
	\label{eq_def_tau}
\end{equation}

  The purpose of this paper is to  develop and  study  a bootstrap suitable
    for  getting approximation of the distribution of the following  class
    of    change-point estimators
\begin{equation}\label{changeest}
    \begin{split}
        \widehat{m}(\gamma)&=\arg\max\{|S_{\gamma}(k)|;k=1,\ldots,n-1\}\\
    &=\min\{k;\,1\ls k<n,\,S_{\gamma}(k)\gs S_{\gamma}(j),\,j=1,\ldots,n-1\},
\end{split}
\end{equation}
where
\[
S_{\gamma}(k)=\left(\frac{n}{k(n-k)}\right)^{\gamma}\sum_{i=1}^k(X(i)-\bar{X}_n)
\]
and $\bar{X}_n=n^{-1}\sum_{i=1}^nX(i)$.

  There is a quite extensive literature concerning  asymptotic behavior of change-point
  estimators for  independent observations.
  For  a survey of  various results, see  e.g.  D\" umb\-gen~\cite{duembgen},
  Cs\" org\H o and Horv\'ath~\cite{csoehor} and Antoch et al.~\cite{anthusest}.
 One of the first papers to derive the limit
distribution
 for $\widehat{m}(1/2)$ and independent errors under local changes has been
  written by Bhattacharya and Brockwell~\cite{bhattbrock76}.
  D\" umbgen~\cite{duembgen} considered  a change in a general model for  AMOC
   with independent observations and developed a  suitable
   bootstrap.
   Antoch et al.~\cite{anthusver95}  studied the asymptotic behavior of
   $S_{\gamma}(k), \, 0\leq \gamma\leq 1/2$, in the model
   (\ref{LPmodel}) with independent identically distributed
    errors and developed and studied a bootstrap valid for
   local changes.  They also obtained various related
results, such as rates of consistency for the estimators and their
limiting distribution. Ferger and Stute~\cite{fergerstute92} and Ferger~\cite{ferger94b,ferger94a}
 studied change-point estimators based on $U$-statistics for  i.i.d. errors.\\

Bai~\cite{bai} and Antoch et al. \cite{anthuslp} analyzed the limit
behavior of various estimators when the error sequence forms  a linear process.
 However, they have not discussed bootstrapping. \\

    Most of the theoretical results concerning  bootstrap methods   in  change-point analysis
   (testing and estimation) have
   been obtained   for  independent observations, see e.g., Hu\v skov\'a~\cite{surveyhus}
.
   Antoch and Hu\v{s}kov\'{a}~\cite{anthusperm} obtained critical
 values for the change-point test related to functionals of $S_{\gamma}(k)$  for
  $H_0:m=n$ ("no change")  vs. $H_1:m<n,d\neq 0$ ("there is a change in the mean") using
permutation methods (or equivalently, bootstrap without
replacement) for the independent case. Recently, Kirch
\cite{kirchphd,kirchblock,kirchfreq}  has   developed  various
bootstrap approximations for critical values  for the above tests of
    "no change" versus  "there is a change in the mean"
   suitable for the  case
   of dependent observations that form a linear process. The results in \cite{kirchblock} can also be modified in a straightforward way for dependent observations as discussed here.

In this paper we develop and prove the validity of a circular overlapping 
block bootstrap for obtaining
asymptotically correct confidence intervals in the case of
dependent errors.

In order to prove validity of the developed bootstrap scheme as well as to obtain the asymptotic under the null hypothesis for the change-point estimator  we
have to use
   some results like laws of (iterated) logarithm or
large numbers for a triangular array. Therefore we need  additionally to the
assumptions (\ref{vor_fehler_1})
 the following one for certain $\delta,\Delta>0$ (in some cases $>2$):
 \begin{enumerate}[($\mathcal{A}$)]\item Let $\{e(i):i\in\mathbb{Z}\}$ be
 a strictly stationary sequence
        with $\E e(0)=0$. Assume there are $\delta,\Delta>0$
        with
        \[
        \E |e(0)|^{2+\delta+\Delta}\ls D_1
        \]
        and 
        \begin{equation}\label{eq_A_2}
            \sum_{k=0}^{\infty}(k+1)^{\delta/2}{\alpha_n}(k)^{\Delta/(2+\delta+\Delta)}
            \ls D_2,
        \end{equation}
        where $\alpha_n(k)$ is the corresponding  strong mixing coefficient,
        i.e.\[ \alpha_n(k)=\sup_{A,B}\left|P(A\cap B)-P(A)P(B)\right|,
        \]
        where $A$ and $B$ vary over the $\sigma$-fields
        $\mathcal{A}(e(0),e(-1),\ldots)$ respectively\\
        $\mathcal{A}(e(k),e(k+1),\ldots)$.
\end{enumerate}

Under this assumption we get moment inequalities (cf. Yokoyama~\cite{yokoyama},
Theorem~1, and Serfling~\cite{serfling}, Lemma B, Theorem 3.1), which in turn
yield laws of large numbers. Moreover the results remain true for triangular arrays that fulfill uniformly the assumptions above. For more details we refer to Kirch~\cite{kirchphd}, Appendix B.2.

In fact we only need this assumption in order to obtain a Donsker type central limit theorem for the partial sums of the errors (to derive the asymptotic under the null hypothesis) as well as bounds on higher order moments of certain sums of the observed error sequence. This in turn yields laws of large numbers and laws of (iterated) logarithm. The proofs can easily be adapted to allow for errors that do not fulfill condition $(\mathcal{A})$ but the necessary moment conditions.

\begin{ex}
Suppose that the errors form a linear process
\[
e(i)=\sum_{s\gs 0}w_s \,\eps(i-s),
\]
where the innovations  $\{\eps(i): -\infty < i <
\infty \}$ are i.i.d. random variables with
\[
	    \E \eps(i) =0,\quad 0<\sigma^2=\E \eps(i)^2<\infty, \quad
	        \E|\eps(i)|^{\nu}<\infty \text{ for some }\nu>4.
	\]
	We suppose that the weights $\{w_s:s\gs 0\}$ satisfy
	\[
		    \sum_{s\gs 0} |w_s|<\infty,\qquad\sum_{s\gs 0} w_s \neq 0.
	    \]
Corollary 4 in Withers~\cite{withers} gives  mild  conditions
under which linear sequences are strong mixing and even provides
the mixing coefficients. This can be used to check 
condition ($\mathcal{A}$). Causal ARMA sequences with
 appropriate innovations, for example, fulfill it for any $\delta,\Delta>0$,
 if the $(2+\delta+\Delta)$-moment of the innovations exists.
\end{ex}

For the sake of simplicity we will only consider the case
$\gamma=1/2$ in the following. The results for $0\ls \gamma <1/2$
can be obtained in a similar way as outlined in Antoch et
al.~\cite{anthusver95}. In the simulation study we will also
consider other choices of $\gamma$, since the asymptotic method
does not give such good approximations for $\gamma<1/2$. The
reason is that the asymptotic distribution in this case depends on
unknown parameters and thus in practice on estimators.

In the present paper we focus on  local alternatives (i.e.,
$d=d_n\to 0$, as $n\to \infty$). To obtain results for
fixed alternatives is more complicated because the limit
distribution of the estimator is determined by finite sums, which
depend on the underlying distribution function. Some comments
concerning the i.i.d. case can be found in the survey paper by
Antoch and \hus~\cite{anthusest}. Furthermore D\"{u}mbgen~\cite{duembgen} 
considers both, local and fixed changes, for independent
observations in a somewhat more general setup, i.e. the parameters
that are subject to change need not be location parameters.

In the following $\wm:=\wm(1/2)$, $S(k):=S_{1/2}(k)$.

\section{Limit Distribution and Rate of Consistency for the Estimators}
\label{sec_asymp}
In this section we summarize and generalize some previous results
by Antoch et al. \cite{anthuslp,anthusver95} that we need in the sequel.

The next theorem gives the rate of consistency for the change-point estimator $\wm$ as well as its limit distribution for a local change.
For the i.i.d. case these results have been obtained by
Antoch et al.~\cite{anthusver95}, Theorem 1 and 2. The second result has been generalized for errors that form a linear process under the additional assumption $\sum_{s\gs 0}\sqrt{s}\,|w_s|<\infty$ by Antoch et al.~\cite{anthuslp} (Theorem~2.2).

\begin{theorem}\label{th_asym_limit}
    Assume that \eqref{LPmodel}-\eqref{eq_def_tau} with $0<\vth<1$ are satisfied and, as $n\to\infty,$
\begin{equation}\label{assump_delta_2}
d_n\to 0,\quad  \left|d_n\right|^{-1}n^{-1/2}(\log n)^{1/2}=o(1).
\end{equation}
Moreover let assumption $(\mathcal{A})$ be fulfilled for some $\delta,\Delta>0$. Then:
\begin{enumerate}[a)]
	\item The following rate of consistency holds, as $n\to\infty$,
		\[
\wm-m=o(d_n^{-2}\log n)\qquad P-a.s.
\]
\item 
    The change-point estimator fulfills, as $n\to\infty$,
    \begin{equation}
        \frac{d_n^2\left( \wm-\lfloor n\vth\rfloor \right)}{\tau^2}\dto \arg\max\{W(t)-|t|/2:t\in \mathbb{R}\},
        \label{eq_asser_limit}
    \end{equation}
    where $\{W(t):t\in\mathbb{R}\}$ is a two-sided Wiener process and $\tau^2$ as in \eqref{eq_def_tau}.
    \end{enumerate}
\end{theorem}
The proof is postponed to Section~\ref{sec_proofs}.

\begin{rem}
We would like to point out that the result in a) also remains true for a fixed change, precisely it suffices to have $\left|d_n\right|^{-1}n^{-1/2}(\log n)^{1/2}=o(1)$ we do not need $d_n\to 0$. As a contrast we do need $d_n\to 0$ to obtain the limit distribution in b), but if this is not fullfilled the proof still shows $\wm-m=O_P(d_n^{-2})$.
\end{rem}

\begin{rem}\label{rem_est_bartlett}Also assertion \eqref{eq_asser_limit} remains true, if one replaces the unknown quantity $\tau$ by a consistent estimator $\widehat{\tau}$. 
	It is also possible to replace $d_n$ by an estimator fulfilling $\widehat{d}_n-d_n=O_P\left(\sqrt{\log n/n}\right)$, since then it follows $\widehat{d}_n/d_n\to 1$ under \eqref{assump_delta_2}.
	\begin{enumerate}[a)]\item The following Bartlett type estimator is a consistent estimator for $\tau^2$ under the conditions of Theorem~\ref{th_asym_limit} (cf. Theorem~1.1 in Berkes et al.~\cite{berkesetal05}, Lemma~A.0.1 in Politis et al.~\cite{politisetal99}, and Theorems~14.1 and 14.5 in Davidson~\cite{davidson02}):
			\begin{equation}\label{eq_est_bartlett}
		 \widetilde{\tau}^2_n(\Lambda_n)=\widehat{R}(0)+2\sum_{k=1}^{\Lambda_n}
		  \left(1-\frac{k}{\Lambda_n}\right) \widehat{R}(k),
	  \end{equation}
	  where
	  \[
	  \widehat{R}(k)=\frac{1}{n}\left(\sum_{t=1}^{\widehat{m}-k}
	  (X_t-\bar{X}_{\widehat{m}})(X_{t+k}-\bar{X}_{\widehat{m}})
	  +\sum_{t=\widehat{m}+1}^{n-k}(X_t-\bar{X}^0_{\widehat{m}})
	  (X_{t+k}-\bar{X}^0_{\widehat{m}} )\right),
	  \]
	  $\widehat{m}=\min\{\arg\max\{\left|S(k)\right|:k=1,\ldots,n\}\}$ and
	  $\bar{X}^0_{\widehat{m}}=\frac{1}{n-\widehat{m}}
	  \sum_{i=\widehat{m}+1}^nX_i$. $\Lambda_n$ should be chosen such that
	  $\Lambda_n^2/\log(n)\to\infty$ and
	  $n^{-1}\Lambda_n^2\log\Lambda_n=o(\log(n))$.

  \item As an estimator for  $d_n$ one can use
	  \begin{equation}\label{eq_est_d}
	  \widehat{d}_n=\frac{1}{n-\wm}\sum_{j=\wm+1}^nX(j)-\frac{1}{\wm}\sum_{j=1}^{\wm}X(j).
	\end{equation}
 Combining Theorem~\ref{th_asym_limit}~a) with Theorem~B.8 and Remark~B.2 in Kirch~\cite{kirchphd} we obtain the rate
 \[
 \widehat{d}_n-d_n=o\left(\sqrt{ \frac{\log n}{n}} \right)\qquad P-a.s.
 \]
 under the assumptions of Theorem~\ref{th_asym_limit}.
	\end{enumerate}
  \end{rem}

\begin{rem}\label{rem_limit_cont}    
	It can be shown that the above limit distribution is continuous
    and explicitly known (confer Remark 2.3 in Antoch et \hus \cite{anthusest}).
Thus the above theorem can be used to construct asymptotic confidence intervals, precisely $(\wm -\widehat{\tau}/\widehat{d}_n q_U(\alpha/2),\wm -\widehat{\tau}/\widehat{d}_n q_L(\alpha/2))$, where $P(U(1/2)\ls q_L(\alpha/2))=\alpha/2$, $P(U(1/2)\gs q_U(\alpha/2))=\alpha/2$ for $U(1/2)=U(\vth,1/2)$ as in Remark~\ref{rem_general_limit} below. Note that $U(\vth,1/2)$ does not depend on the unknown parameter $\vth$.
\end{rem}

\begin{rem}\label{rem_general_limit}
For $0\ls \gamma <\frac 1 2$ the limit distribution depends on the unknown parameter~$\vth$. Precisely it can be shown that the limit is
\begin{align*}
&U(\vth,\gamma)=\arg\max\{W(t)-|t|g_{\vth,\gamma}(t)\},\\
&\text{where}\quad g_{\vth,\gamma}(t)=\begin{cases} (1-\vth)(1-\gamma)+\vth\gamma,&t<0,\\
	(1-\vth)\gamma+\vth(1-\gamma),&t\gs 0.\\
\end{cases}
\end{align*}
Remark 2.3 in Antoch et \hus \cite{anthusest} gives a closed  formula for the limit distribution. We would like to point out a small (but for simulations very important) misprint there: The integral is equal to $c_2/(c_1+c_2)+\ldots$ and not to $c_1/(c_1+c_2)+\ldots$.
\end{rem}

\section{Bootstrap approximations}

\label{sec_mainresults}
Antoch et al.~\cite{anthusver95}  propose a bootstrap with replacement
of the estimated error sequence to obtain confidence intervals for the
change-point. Since in our case the error sequence is no longer independent
 we  have to use a slightly different approach here. We still bootstrap
 the estimated error sequence with replacement, but we will now use a
  circular moving block bootstrap as suggested by Politis and
  Romano~\cite{politisromano1992}. It has the advantage over the
  regular moving blocks bootstrap by K\"{u}nsch~\cite{kuensch}
  that the sample mean is unbiased. Another possibility is to
   use non-overlapping blocks as suggested by Carlstein~\cite{carlstein},
   but this bootstrap does behave slightly worse in simulations.

Kirch~\cite{kirchphd,kirchblock} used block bootstrapping procedures (more precisely a block permutation method as well as a circular and non-circular block bootstrap) to get approximations for the critical values of the change-point test corresponding to the above problem.

Block bootstrapping methods split the observation sequence of length~$n$ into
sequences of length~$K$. Then we put $L$ of them together to a bootstrap sequence (i.e. $n=K L$). We keep the order within the blocks. $K$ and $L$ depend on $n$ and
converge to infinity with $n$.

The idea is that, for properly chosen block-length $K$, the block contains enough information about the
dependency structure so that the estimate is close to the null
hypothesis.

We assume in the following that
\begin{equation}\label{assum_KL}K,L\to\infty\quad\text{ and}\quad K=K(L),\;n=n(L)=K
L,\; K/L=o(1).
\end{equation}

Let $\wmu_1$ be an estimator for $\mu$,
$\wmu_2$ for $\mu+d_n$, and $\widehat{d}_n$ for $d_n$, e.g.
\begin{align}\label{eq_wmu}
&{\wmu}_1=\frac{1}{\wm}\sum_{i=1}^{\wm}X(i),\qquad
{\wmu}_2=\frac{1}{n-\wm}\sum_{i=\wm+1}^{n}X(i),\qquad \widehat{d}_n=\wmu_2-\wmu_1,
\end{align}
where $\wm=\wm(1/2)$ as in \eqref{changeest}.
Remark~\ref{rem_est_bartlett} yields that $\widehat{d}_n$  fulfills assumption~\eqref{assump_delta_est_2}.

Define the estimated
residuals and the centered residuals by
\[
\widehat{e}(i)=X(i)-\wmu_1 1_{\{ i\ls \wm\}}-\wmu_2 1_{\{\wm< i\ls n\}},
\]
\[
\we(i)=\widehat{e}(i)-\frac 1 n \sum_{j=1}^n\widehat{e}(j),
\]
respectively. Throughout the paper the following representation
will turn out to be very useful
\begin{equation}
    \begin{split}
        \we(i)=&e(i)-\bar{e}_n+d_n\left( 1_{\{m<i\ls n\}}-\frac{n-m}{n} \right)-\widehat{d}_n\left( 1_{\{\wm<i\ls n\}}-\frac{n-\wm}{n} \right).
    \end{split}
    \label{eq_decomp_etilde}
\end{equation}

Let $(U(1),\ldots,U(L))$ be i.i.d. with $P(U(1)=i)=\frac 1 n$ for
$i=0,\ldots,n-1$ independent of the observations $X(1),\ldots,X(n)$. Take the i.i.d. bootstrap sample
$e^*(Kl+k)=\we(U(l)+k)$, where $\we(j)=\we(j-n)$ for $j>n$ (hence the name circular bootstrap).

Consider the bootstrap observations
\[
X^*(Kl+k)=e^*(Kl+k)+\wmu_1 1_{\{Kl+k\ls \wm\}}+\wmu_2 1_{\{\wm < Kl+k\ls n\}}.
\]
We  now deal with the following bootstrap estimator of the change-point $m$
\begin{equation}\label{eq_boot_est}
\wm^*=\arg\max\{|S^*(k)|;k=1,\ldots,n-1\},
\end{equation} where
\[
S^*(k)=\left(\frac{n}{k(n-k)}\right)^{1/2}\sum_{i=1}^k(X^*(i)-\bar{X}^*_n).
\]

Now we are ready to present results on the asymptotic
 behavior of the bootstrap estimator $\wm^*$ defined in
  \eqref{eq_boot_est} of the change-point together with a
  short discussion, how to apply the result to obtain confidence intervals.

With $P^*$, $\E^*$, $\var^*,\ldots$ we will denote probability,
expectation, variance,$\ldots$, given\\ $X(1),\ldots,X(n)$.

\begin{theorem}\label{theorem_local}
    Assume that \eqref{LPmodel}-\eqref{eq_def_tau} with $0<\vth<1$ and \eqref{assum_KL} hold. Moreover let
\begin{equation}
    \label{assump_delta_est_2}
    \widehat{d}_n-d_n=O\left(\sqrt{\frac{\log n}{n}}\right)\qquad P-a.s.
\end{equation}
be fulfilled in addition to \eqref{assump_delta_2}.
    Moreover let assumption $(\mathcal{A})$ be fulfilled for some
    \begin{equation}
        \label{assum_deltas_main}
        0<\delta^{(2)}+\Delta^{(2)}<(\delta^{(1)}-2)/2<(\nu-4)/2$, $\Delta^{(1)}=\nu-2-\delta^{(1)}.
    \end{equation}
If $d_n^2 K\to 0$, then
\begin{align*}
&\sup_{x\in\mathbb{R}}\left|P^*(\widehat{d}^2_n(\wm^*-\wm)\ls
x)-P(d^2_n(\wm-m)\ls x)\right|\to 0\qquad P-a.s.
\end{align*}
\end{theorem}

Since the limit distribution (for both the bootstrap as well as null asymptotic)
 is continuous (as has been pointed out by Remark~\ref{rem_limit_cont})
  the described  sampling scheme provides  bootstrap approximations
  to the $(1-\alpha)$-quantile for arbitrary $\alpha\in (0,1)$.
Thus the bootstrap based approximation for the change-point $m$ can be constructed
 along the usual lines.
Precisely the $(1-\alpha)$-bootstrap confidence interval is given by
\[
[2\wm-q_L^*(\alpha/2),2\wm-q_U^*(\alpha/2)],
\]
where
$$
q_L^*(\alpha/2)=\sup\{ u;\, P^*\left( \wm^* < u\right)\ls\alpha/2\}
$$
and
$$
q_U^*(\alpha/2)=\inf\{ u;\, P^*\left( \wm^* >
u\right)\ls\alpha/2\}.
$$
Usually one uses the empirical bootstrap distribution of $\wm^*$ for say $10\,000$ random bootstrap samples. Further discussions on
bootstrap approximations of confidence intervals (for the similar
case of i.i.d. errors) can be found in Antoch and \hus
\cite{anthusest}.

\begin{rem}
    There are also several other possibilities of bootstrapping.
     For example we can use a non-circular approach and/or
     non-overlapping blocks. Simulations for the bootstrap where $(U(1),\ldots,U(L))$
      are i.i.d. uniformly distributed on $\{0,\ldots,L-1\}$ and $e^*(Kl+k)=\we(K U(l)+k)$
      indicate that this bootstrap does not perform quite as good as the bootstrap proposed above.

\end{rem}

\section{Simulation Study}

In the previous chapter we have established the asymptotic validity of the bootstrap confidence intervals. The question remains how well these confidence intervals behave for small samples and also how well they behave in comparison with the asymptotic intervals. 

In this section we  not only consider $\gamma=\frac 1 2$ but also $\gamma=0$. The important difference is that the asymptotic confidence intervals depend on the unknown parameter $\vth$ for $0\ls \gamma<\frac 1 2$. Not surprisingly it  turns out that the asymptotic intervals behave  better for $\gamma=\frac 1 2$, whereas in all other cases it is better to use the bootstrap intervals.\\
Moreover we  consider changes in the mean of $d= 0.5,1,2,4$. The latter ones can hardly be regarded as local changes, however we are still interested in the behavior of the bootstrap intervals, since we conjecture it will also be valid in those cases.

For the simulations we use an autoregressive sequence of order one as an error sequence with standard normally distributed innovations and different values of $\rho=0.1,0.3$. We consider changes at $m=20,40$. We use the estimator $\widehat{d}_n$ as in \eqref{eq_est_d} respectively \eqref{eq_wmu}
 and - for the asymptotic method - the Bartlett estimator given in \eqref{eq_est_bartlett} with $\Lambda_n=0.1\,n$, because in the simulation study conducted by Antoch et al.~\cite{anthuslp} this choice gave best results in the AR(1)-case.

The goodness of confidence intervals can essentially be determined by two criteria:
\begin{enumerate}[$\mathcal{C}$.1]
	\item The probability that the actual change-point is outside the (1-$\alpha$)-confidence interval should be close to (smaller than) $\alpha$.
	\item The confidence intervals should be short.
\end{enumerate}
We visualize the first quantity by using CoLe-Plots ({\bf Co}nfidence-{\bf Le}vel-Plots) and the second one by using CoIL-Plots ({\bf Co}nfidence-{\bf I}nterval-{\bf L}ength-Plots).

In fact we have done more simulations (such as QQ-plots or tables of the quantiles of $m$) for a large amount of different combinations of  parameters as well as different possible bootstrap procedures. 
The problem, however, is that for the bootstrap they  only give result for one specific underlying sequence and are thus rather not as informative.  For this reason and also due to similarity of results as well as due to limitations of space we restrict ourselves to the following plots.\\[2mm]
{\bf CoLe-Plots}\\[2mm]
We  explain how the plots are created using the example of asymptotic confidence intervals. The general version of Theorem~\ref{th_asym_limit}
yields that the asymptotic confidence intervals are  calculated using the distribution of $Z=\widehat{m}-{\widehat{\tau}^2}/{\widehat{d}^2_n}V(\widehat{\vth},\gamma)$, where $V(\widehat{\vth},\gamma)$ is as in Remark~\ref{rem_general_limit} and $\widehat{\vth}=\wm/n$. 
Note that
\begin{align*}
m\not\in CI(1-\alpha)&\quad\Longleftrightarrow \quad P(Z\ls m)\ls \alpha/2\quad\text{or}\quad P(Z\gs m)\ls \alpha/2\\
&\quad\Longleftrightarrow \quad 2\min(P(Z\ls m), P(Z\gs m))\ls \alpha.
\end{align*}
The CoLe-Plots now draw the empirical distributions function (based on $1\,000$ observation sequences) of $2\min(P(Z\ls m),P(Z\gs m))$.
\\
Thus for given $\alpha$ on the $x$-axis the plot shows the empirical probability that $m$ is outside the $(1-\alpha)$-confidence interval on the $y$-axis, hence it visualizes $\mathcal{C}.1$. Optimally, the plot should be below or (even better) on the diagonal.

For the bootstrap confidence intervals the procedures works exactly the same but now the intervals are calculated using the (empirical, based on $10\,000$ resamples) distribution of $\widetilde{Z}=2\widehat{m}-\widehat{m}^*$.
\\[2mm]
{\bf CoIL-Plots}\\[2mm]
We calculate for $1\,000$ observation sequences the length of the confidence intervals for levels $\alpha=0.01,0.02,\ldots, 0.1$. The empirical bootstrap distribution is based on $10\,000$ random samples as before. Then we plot the mean using a  thick line (as well as the upper and lower quartiles with  thin lines), linearly interpolated. So these plots visualize the length of the intervals and thus  $\mathcal{C}.2$.
\\[2mm]
Note that the scale on the $y$-axis is not the same for different pictures. This way we can better compare the asymptotic with the bootstrap method.

\sublabon{figure}

\begin{figure}
\subfigure[CoLe-Plot: $d=0.5$]{
\begin{minipage}[b]{.46\linewidth}
\begin{turn}{270}
\includegraphics[width=4.8cm]{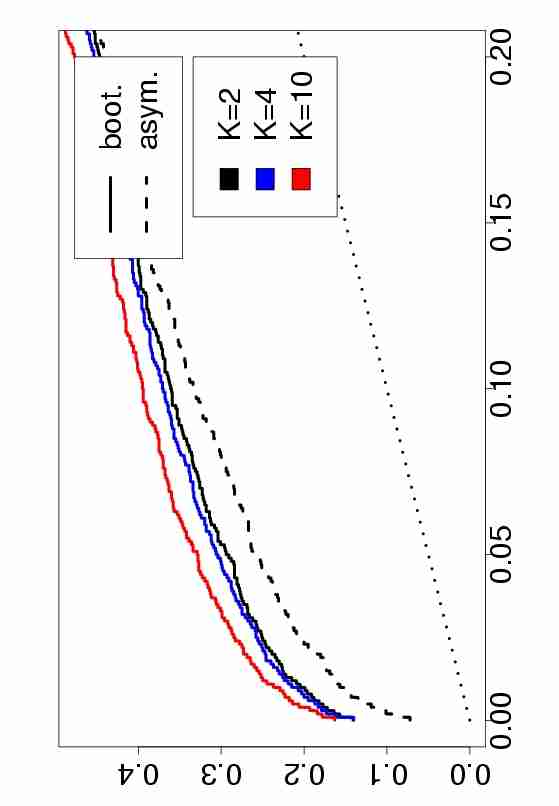}
\end{turn}
\end{minipage} \hfill}\quad
\subfigure[CoIL-Plot: $d=0.5$]{
\begin{minipage}[b]{.46\linewidth}
\begin{turn}{270}
\includegraphics[width=4.8cm]{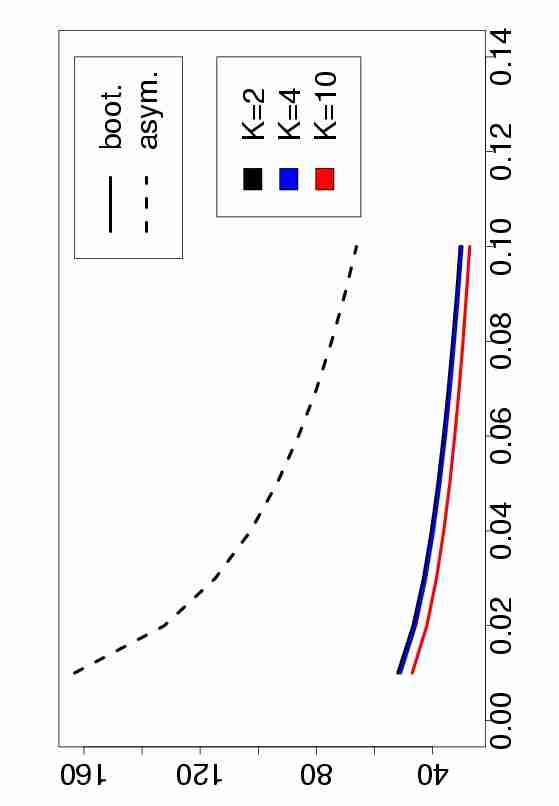}
\end{turn}
\end{minipage} \hfill}\quad
\subfigure[CoLe-Plot: $d=1$]{
\begin{minipage}[b]{.46\linewidth}
\begin{turn}{270}
\includegraphics[width=4.8cm]{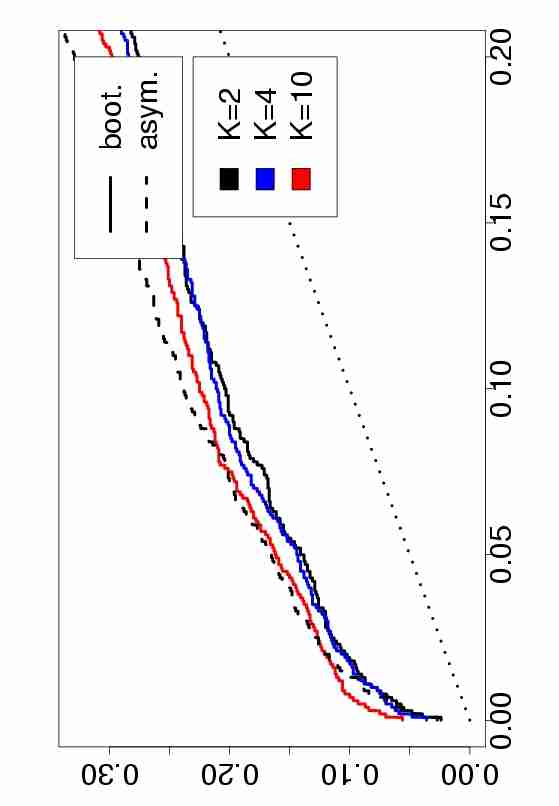}
\end{turn}
\end{minipage} \hfill}\quad
\subfigure[CoIL-Plot: $d=1$]{
\begin{minipage}[b]{.46\linewidth}
\begin{turn}{270}
\includegraphics[width=4.8cm]{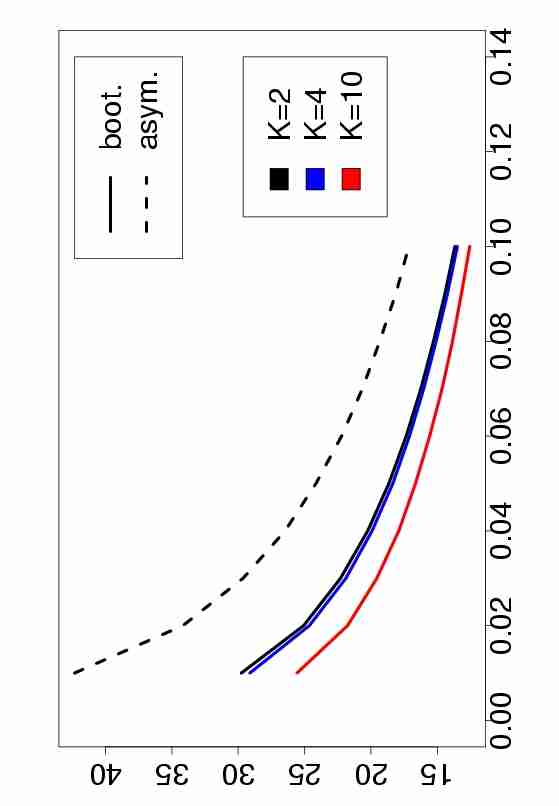}
\end{turn}
\end{minipage} \hfill}\quad
\subfigure[CoLe-Plot: $d=2$]{
\begin{minipage}[b]{.46\linewidth}
\begin{turn}{270}
\includegraphics[width=4.8cm]{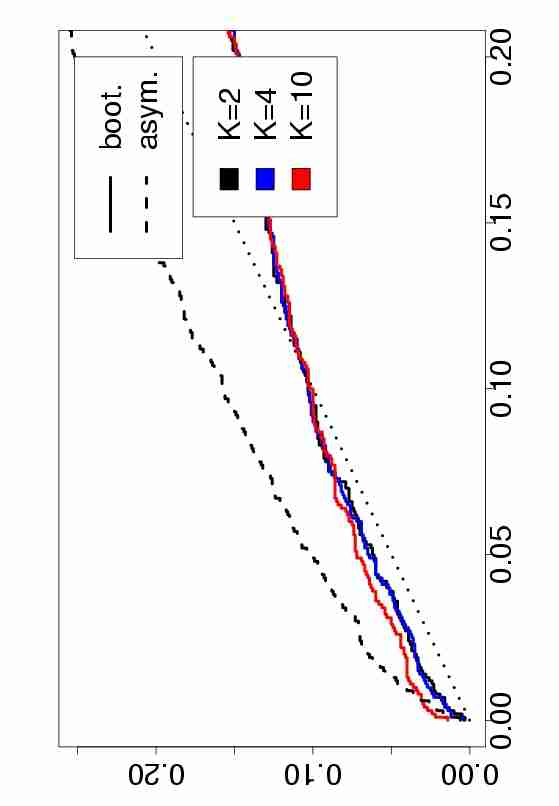}
\end{turn}
\end{minipage} \hfill}\quad
\subfigure[CoIL-Plot: $d=2$]{
\begin{minipage}[b]{.46\linewidth}
\begin{turn}{270}
\includegraphics[width=4.8cm]{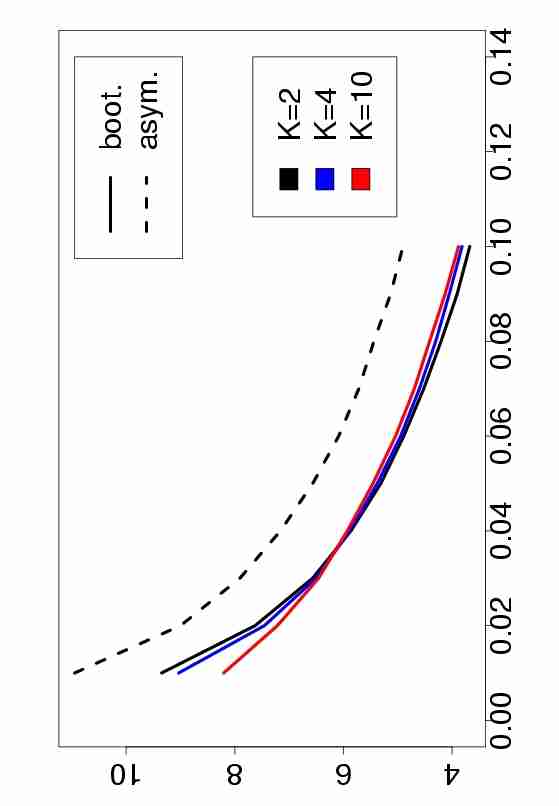}
\end{turn}
\end{minipage} \hfill}\quad
\subfigure[CoLe-Plot: $d=4$]{
\begin{minipage}[b]{.46\linewidth}
\begin{turn}{270}
\includegraphics[width=4.8cm]{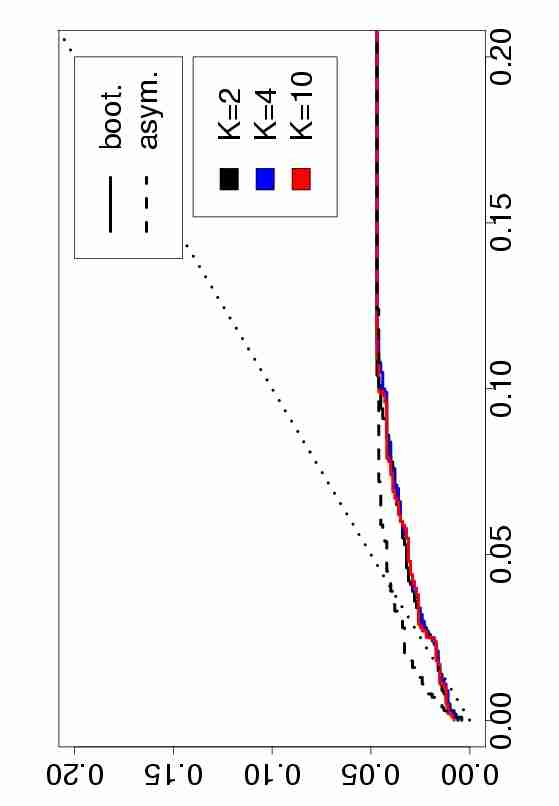}
\end{turn}
\end{minipage} \hfill}\quad
\subfigure[CoIL-Plot: $d=4$]{
\begin{minipage}[b]{.46\linewidth}
\begin{turn}{270}
\includegraphics[width=4.8cm]{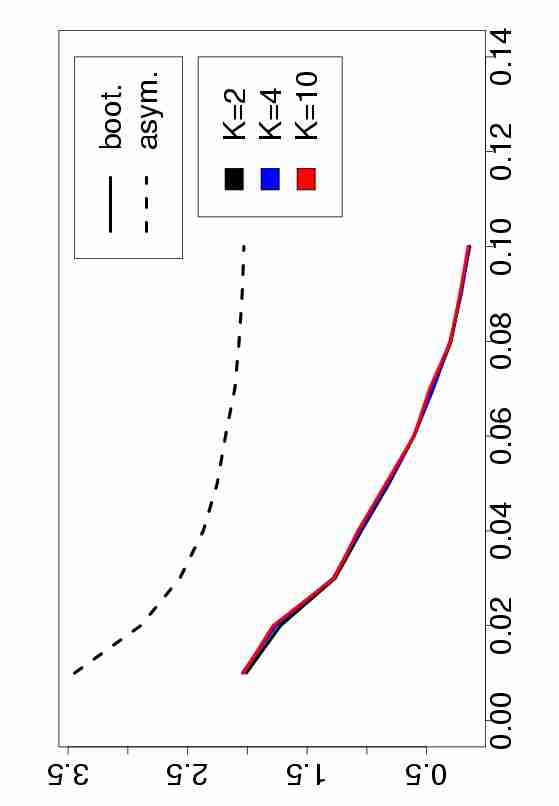}
\end{turn}
\end{minipage} \hfill}\quad
\caption{CoLe- and CoIL-Plots for $n=80$, $m=40$, $\rho=0.1$, and $\gamma=0$ for the asymptotic method as well as bootstrap method with different block-lengths}\label{figure_spcplot4a}
\end{figure}

\begin{figure}
\subfigure[CoLe-Plot: $d=0.5$, $m=40$]{
\begin{minipage}[b]{.46\linewidth}
\begin{turn}{270}
\includegraphics[width=4.8cm]{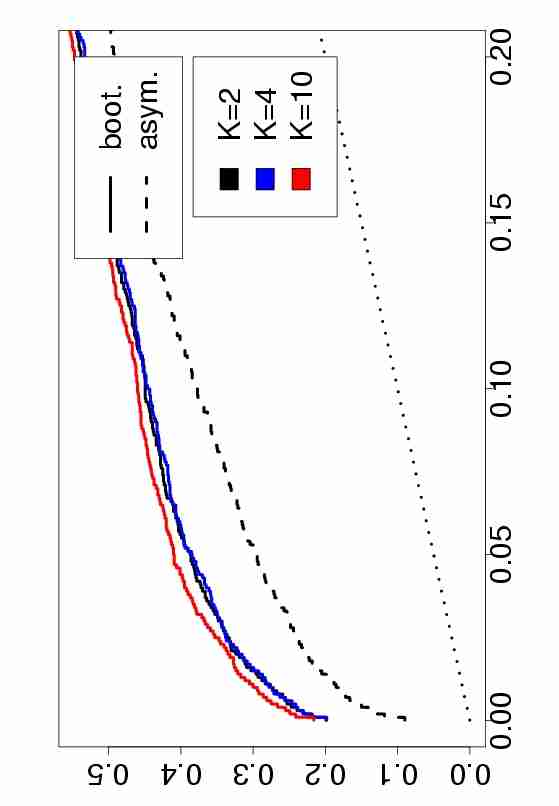}
\end{turn}
\end{minipage} \hfill}\quad
\subfigure[CoIL-Plot: $d=0.5$, $m=40$]{
\begin{minipage}[b]{.46\linewidth}
\begin{turn}{270}
\includegraphics[width=4.8cm]{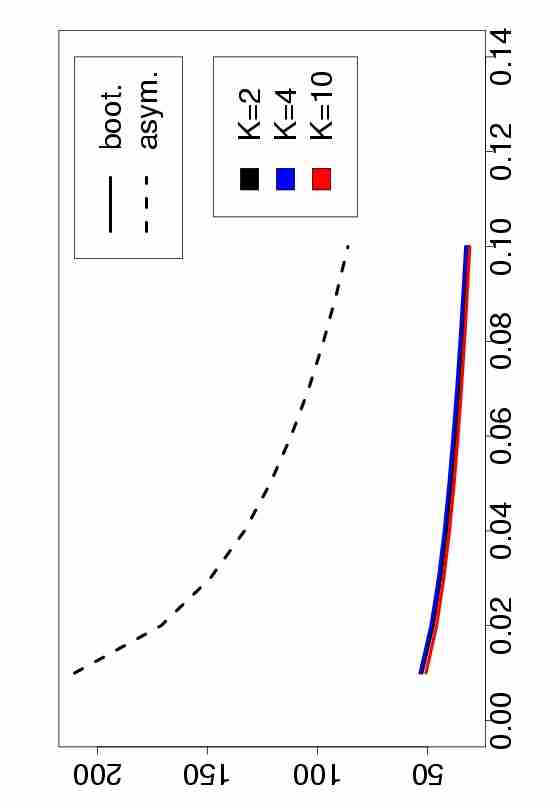}
\end{turn}
\end{minipage} \hfill}\quad
\subfigure[CoLe-Plot: $d=2$, $m=40$]{
\begin{minipage}[b]{.46\linewidth}
\begin{turn}{270}
\includegraphics[width=4.8cm]{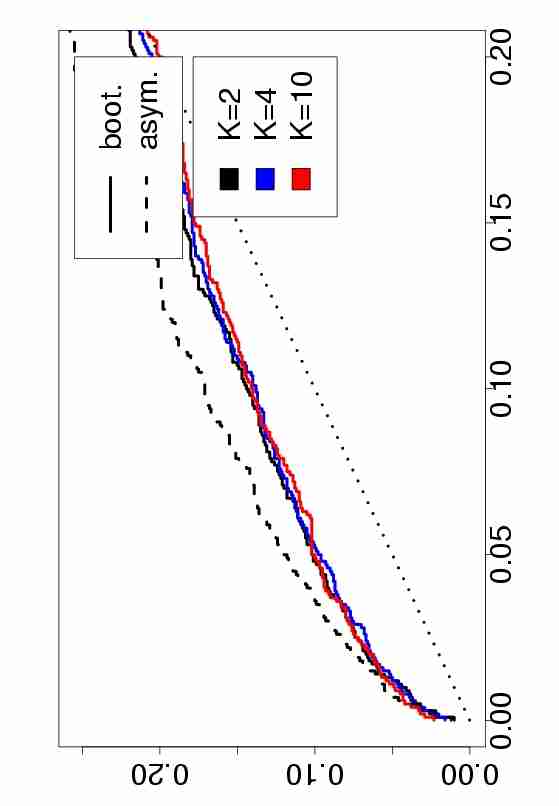}
\end{turn}
\end{minipage} \hfill}\quad
\subfigure[CoIL-Plot: $d=2$, $m=40$]{
\begin{minipage}[b]{.46\linewidth}
\begin{turn}{270}
\includegraphics[width=4.8cm]{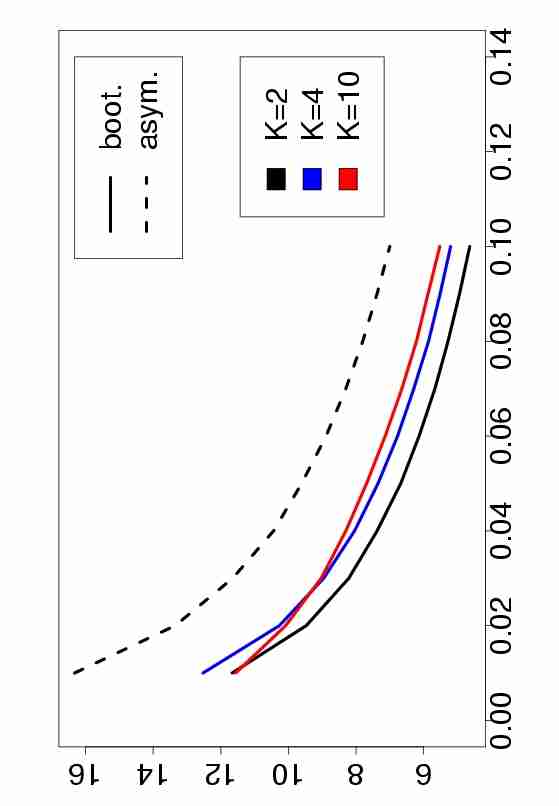}
\end{turn}
\end{minipage} \hfill}\quad
\subfigure[CoLe-Plot: $d=0.5$, $m=20$]{
\begin{minipage}[b]{.46\linewidth}
\begin{turn}{270}
\includegraphics[width=4.8cm]{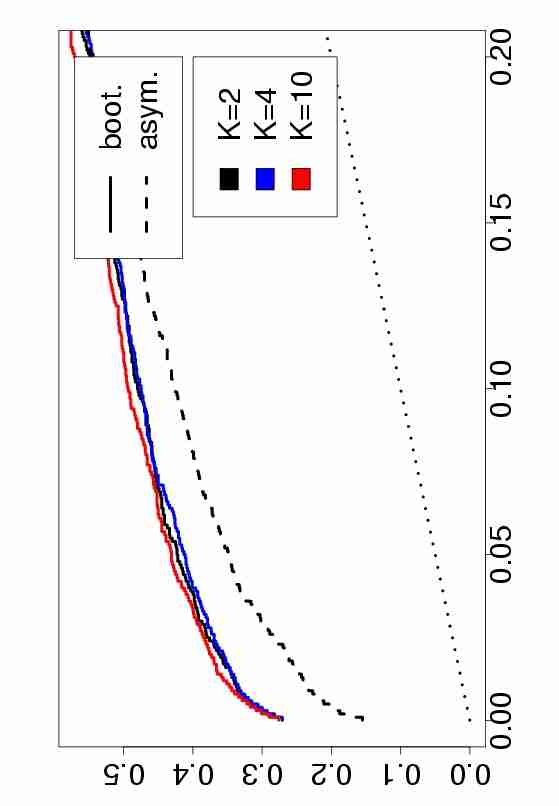}
\end{turn}
\end{minipage} \hfill}\quad
\subfigure[CoIL-Plot: $d=0.5$, $m=20$]{
\begin{minipage}[b]{.46\linewidth}
\begin{turn}{270}
\includegraphics[width=4.8cm]{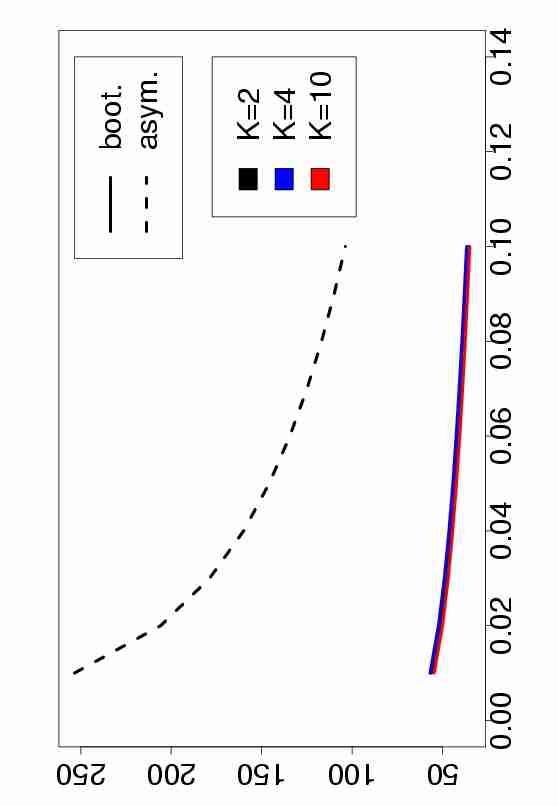}
\end{turn}
\end{minipage} \hfill}\quad
\subfigure[CoLe-Plot: $d=2$, $m=20$]{
\begin{minipage}[b]{.46\linewidth}
\begin{turn}{270}
\includegraphics[width=4.8cm]{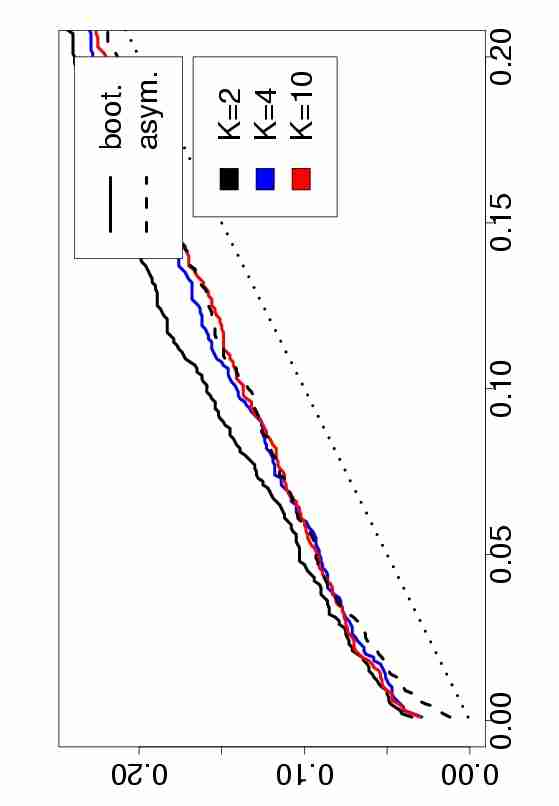}
\end{turn}
\end{minipage} \hfill}\quad
\subfigure[CoIL-Plot: $d=2$, $m=20$]{
\begin{minipage}[b]{.46\linewidth}
\begin{turn}{270}
\includegraphics[width=4.8cm]{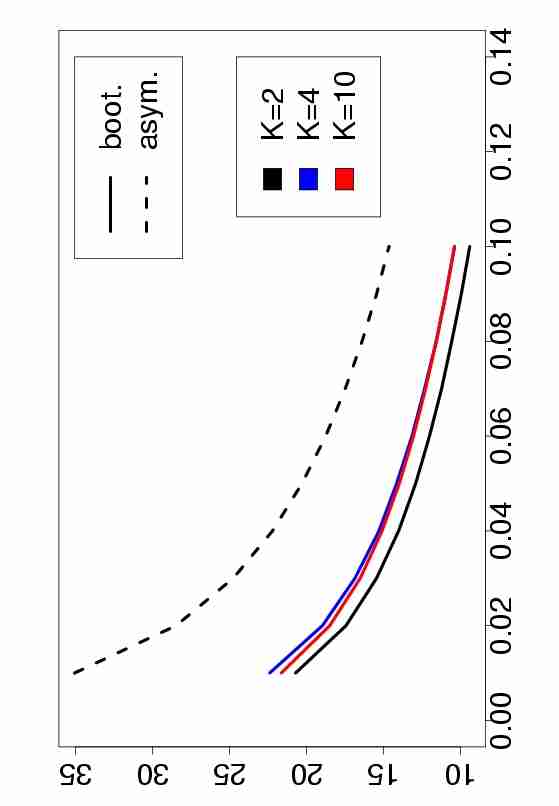}
\end{turn}
\end{minipage} \hfill}\quad
\caption{CoLe- and CoIL-Plots for $n=80$, $\rho=0.3$, and $\gamma=0$ for the asymptotic method as well as bootstrap method with different block-lengths}
\label{figure_spcplot4d}\end{figure}

\sublaboff{figure}

\begin{figure}
\subfigure[CoLe-Plot: $d=0.5$]{
\begin{minipage}[b]{.46\linewidth}
\begin{turn}{270}
\includegraphics[width=4.8cm]{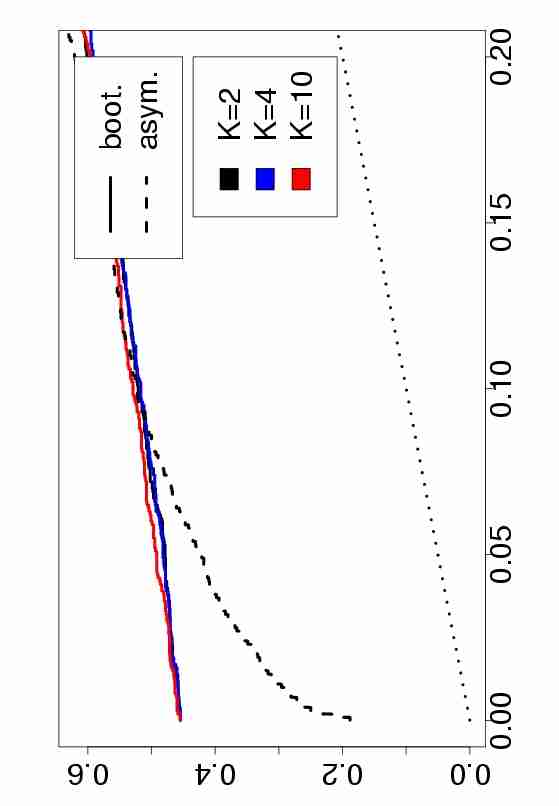}
\end{turn}
\end{minipage} \hfill}\quad
\subfigure[CoIL-Plot: $d=0.5$]{
\begin{minipage}[b]{.46\linewidth}
\begin{turn}{270}
\includegraphics[width=4.8cm]{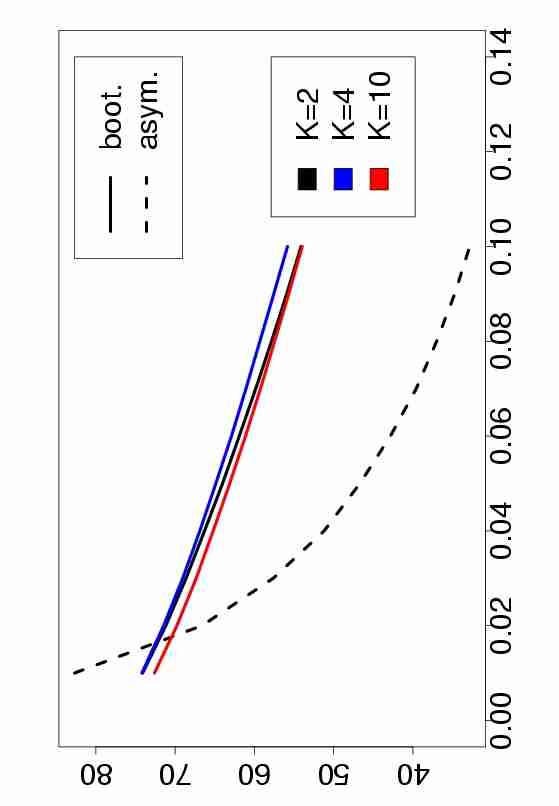}
\end{turn}
\end{minipage} \hfill}\quad
\subfigure[CoLe-Plot: $d=1$]{
\begin{minipage}[b]{.46\linewidth}
\begin{turn}{270}
\includegraphics[width=4.8cm]{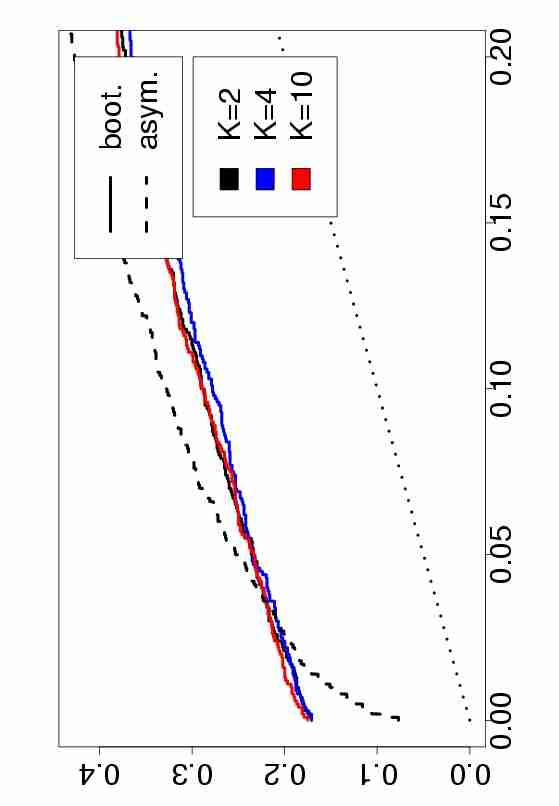}
\end{turn}
\end{minipage} \hfill}\quad
\subfigure[CoIL-Plot: $d=1$]{
\begin{minipage}[b]{.46\linewidth}
\begin{turn}{270}
\includegraphics[width=4.8cm]{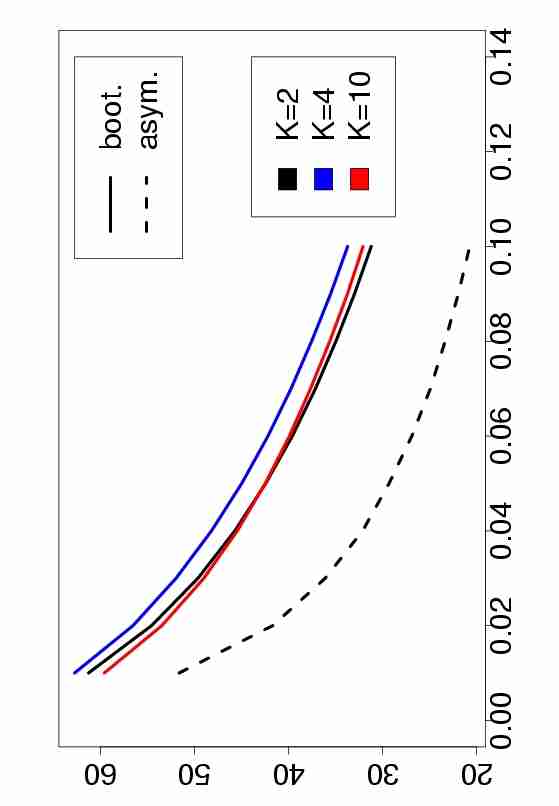}
\end{turn}
\end{minipage} \hfill}\quad
\subfigure[CoLe-Plot: $d=2$]{
\begin{minipage}[b]{.46\linewidth}
\begin{turn}{270}
\includegraphics[width=4.8cm]{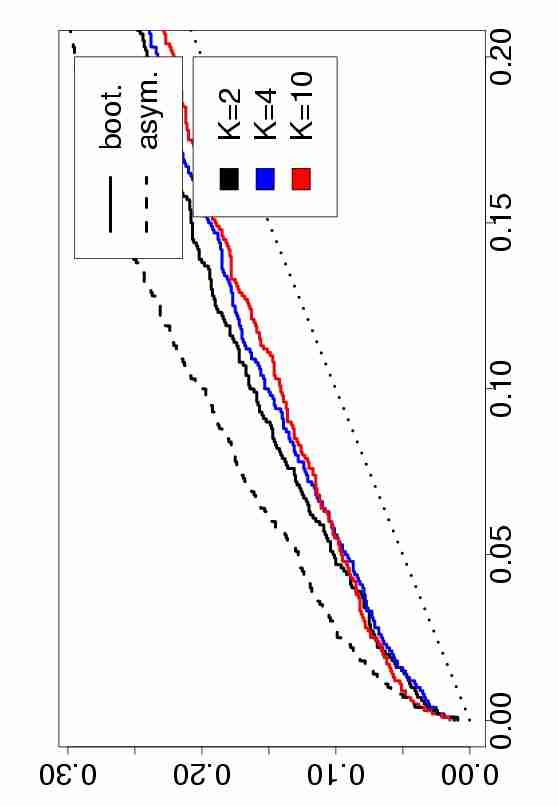}
\end{turn}
\end{minipage} \hfill}\quad
\subfigure[CoIL-Plot: $d=2$]{
\begin{minipage}[b]{.46\linewidth}
\begin{turn}{270}
\includegraphics[width=4.8cm]{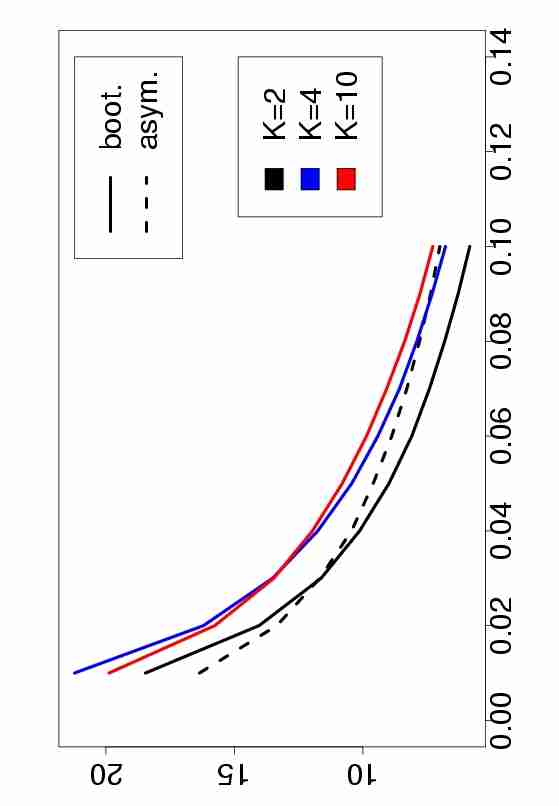}
\end{turn}
\end{minipage} \hfill}\quad
\subfigure[CoLe-Plot: $d=4$]{
\begin{minipage}[b]{.46\linewidth}
\begin{turn}{270}
\includegraphics[width=4.8cm]{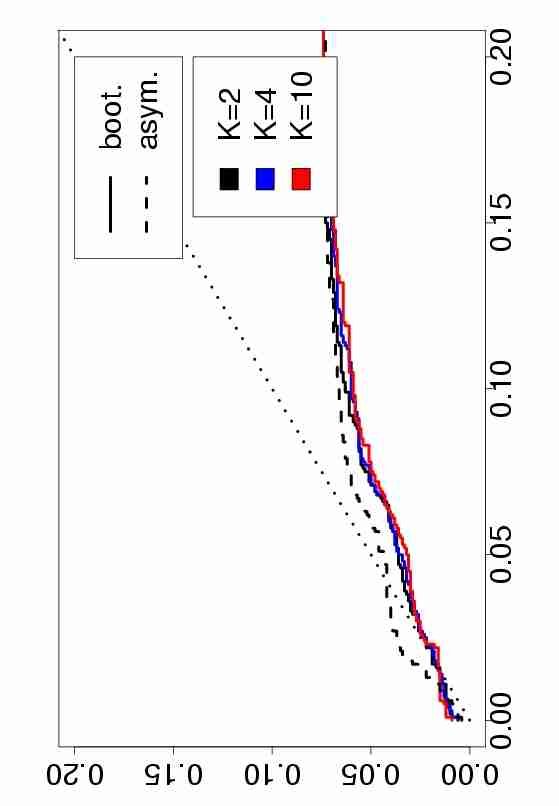}
\end{turn}
\end{minipage} \hfill}\quad
\subfigure[CoIL-Plot: $d=4$]{
\begin{minipage}[b]{.46\linewidth}
\begin{turn}{270}
\includegraphics[width=4.8cm]{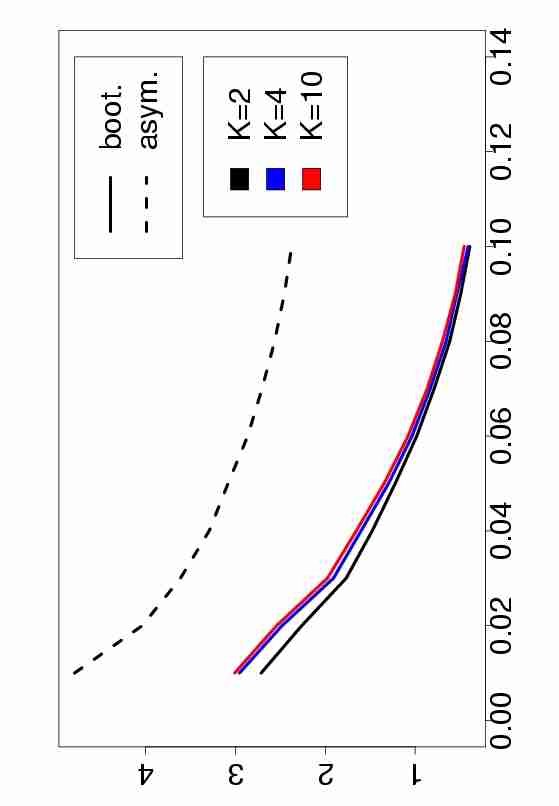}
\end{turn}
\end{minipage} \hfill}\quad
\caption{CoLe- and CoIL-Plots for $n=80$, $m=20$, $\rho=0.3$, and $\gamma=0.5$ for the asymptotic method as well as bootstrap method with different block-lengths}\label{figure_spcplot4h}
\end{figure}
\begin{figure}\centering
\begin{turn}{270}
\includegraphics[width=10cm]{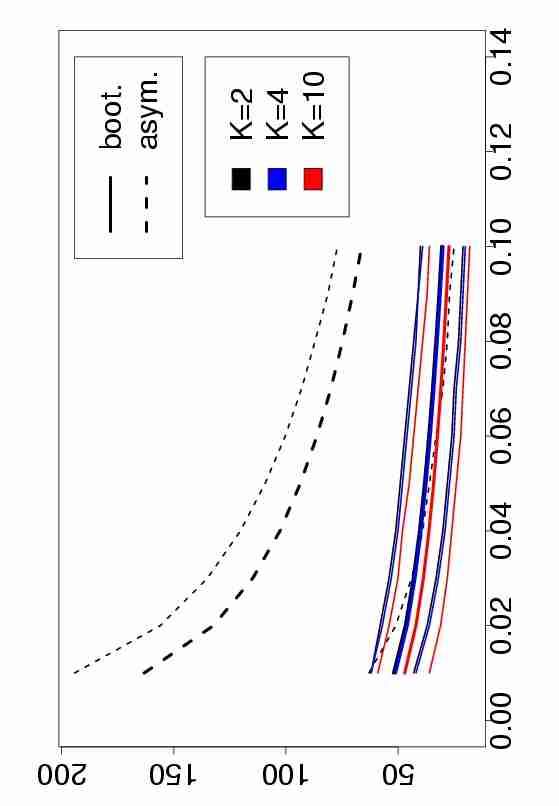}
\end{turn}
\caption{CoIL-Plot for $n=80$, $m=40$, $\rho=0.1$, $\gamma=0$, and $d=0.5$, which additionally to Figure~\ref{figure_spcplot4a}~(2) includes quartiles of the interval length}
	\label{figure_CoIL}
\end{figure}

The plots are given in figures \ref{figure_spcplot4a}-\ref{figure_spcplot4d} for $\gamma=0$ and \ref{figure_spcplot4h} for $\gamma=0.5$. Concerning the CoIL-Plots we only plot the means for better readability. In Figure~\ref{figure_CoIL} we give the CoIL-Plot corresponding to Figure~\ref{figure_spcplot4a}~(2) including the quartiles to give a better idea of the distribution of the  length of the confidence interval.
Concerning $\gamma=0$ we see that for small $d_n$ the actual cover probability of the interval is too small for both methods, yet the asymptotic interval is somewhat better than the bootstrap intervals. At the same time  the length of the asymptotic interval is very large, much larger than the length of the bootstrap interval. Frequently it is even longer than the observation sequence. We did not correct upper and lower bounds of the intervals by $0$ respectively $80$, but bootstrap intervals can also be outside that possible range. 

In fact it is somewhat surprising that even though the intervals are quite long the levels are not as good. The reason is that the change-point estimator for such a small change (and relatively few observations points) is frequently not very good. A typical example is an observation sequence with a change at $m=20$, where the estimator suggests a change at $70$. This results in intervals that do not contain the actual change-point. Also this leads to a wrong estimation of the parameters of the underlying asymptotic distribution, which is then highly skewed in the wrong direction. Thus the lower quantile of the interval is something around $60$, whereas the upper quantile is far bigger than $80$.

For more obvious changes the level of the intervals as well as the length becomes better. This is somewhat surprising in case of the asymptotic intervals  because for fixed changes the asymptotic is not valid. The reason is that we  have an interval around the change-point estimator, which is quite good for more obvious changes.

If the changes are closer to the border of the interval, the levels for both methods deteriorate somewhat. The same holds true for stronger correlation of the underlying error sequence.

Overall the bootstrap intervals behave better than the asymptotic intervals.

However, in the case of $\gamma=0.5$ the asymptotic distribution does not depend on unknown parameters anymore. In this case the asymptotic confidence intervals for local changes are in fact better than the bootstrap intervals. The levels of both methods are for small $d$ somewhat worse than for $\gamma=0$, but the lengths are much better, especially for the asymptotic intervals. However, for more obvious changes the bootstrap intervals are again better than the asymptotic ones. This is due to the fact that the asymptotic does not hold in this case.

It is worth noting that the performance of the bootstrap method does not seem to depend significantly on the choice of the block-length. This is in contrast to the situation where we bootstrap critical values for change-point tests (cf. Kirch~\cite{kirchblock}) where a larger block-length was needed when the data was more dependent.

In real-life situations we recommend to rather use the bootstrap intervals, since they work no matter what $\gamma$ and for both, local as well as fixed changes.

\section{Proofs}\label{sec_proofs}
Throughout the proofs we use the notation $a_n\ll b_n$ for $a_n=O(b_n)$.

We start with the proof of Theorem~\ref{th_asym_limit} in Section~\ref{sec_asymp}.

\begin{proof}[of Theorem~\ref{th_asym_limit}]
    We only sketch the proof, because it is very similar to the proof of Theorem 1 respectively 2 in Antoch et al.~\cite{anthusver95}. First note that
    \[
          \wm=\arg\max\{V(k):k=1,\ldots,n-1\},\qquad \text{where }V(k)=S(k)^2-S(m)^2.
    \]
Simple calculations yield for $k<m$
\begin{equation}\label{eq_decomp_V}
    \begin{split}
    V(k)&= \frac{n(m-k)(n-m-k)}{k(n-k)m(n-m)}\left( \sum_{i=1}^k(e(i)-\bar{e}_n) \right)^2\\
    &\;-\frac{n}{m(n-m)}\sum_{i=k+1}^m(e(i)-\bar{e}_n)\left( \sum_{i=1}^k(e(i)-\bar{e}_n)+\sum_{i=1}^m(e(i)-\bar{e}_n) \right)\\
    &\;+2d_n \sum_{i=k+1}^m(e(i)-\bar{e}_n)\\
    &\;+2d_n \frac{m-k}{n-k}\sum_{i=1}^k(e(i)-\bar{e}_n)\\
    &\;+d_n^2\frac{(n-m)(k-m)}{n-k}\\
    &= A_{k1}+A_{k2}+A_{k3}+A_{k4}+A_{k5}.
\end{split}
\end{equation}
First we show assertion~a), i.e. the $a.s.$ rate of consistency for the change-point estimator.
Theorem B.8 b) and Remark B.2 in Kirch~\cite{kirchphd} give
\begin{align*}
    &\max_{k\ls m}\left|\frac{A_{k1}}{A_{k5}}\right|=\max_{k\ls m-a_n}d_n^{-2}\frac{n|n-m-k|}{m(n-m)^2}\left| \frac{1}{\sqrt{k}}\sum_{i=1}^k(e(i)-\bar{e}_n) \right|^2\\
    &=o\left(d_n^{-2}n^{-1}\log n  \right)=o(1)\qquad P-a.s.
\end{align*}
Similarly we get for $a_n=\eps d_n^{-2}\log n$, where $\eps>0$ is an arbitrary fixed constant,
{\allowdisplaybreaks
\begin{align*}
    &\max_{k\ls m-a_n}\left|\frac{A_{k2}}{A_{k5}}\right|=\max_{k\ls m-a_n}d_n^{-2}\frac{n(n-k)}{m(n-m)^2}\left| \frac{1}{{m-k}}\sum_{i=k+1}^m(e(i)-\bar{e}_n)\right|\\*
    &\phantom{ \max_{k\ls m-a_n}\left|\frac{A_{k2}}{A_{k5}}\right|=\max_{k\ls m-a_n} }\quad\cdot\left| \sum_{i=1}^k(e(i)-\bar{e}_n)+ \sum_{i=1}^m(e(i)-\bar{e}_n)\right| \\*
    &\qquad=o\left(d_n^{-2}\sqrt{\frac{\log n}{n}}\sqrt{\frac{\log n}{a_n}} \right)=o(1)\qquad P-a.s.,\\
    &\max_{k\ls m-a_n}\left|\frac{A_{k3}}{A_{k5}}\right|=\max_{k\ls m-a_n}2d_n^{-1}\frac{n-k}{n-m}\left| \frac{1}{{m-k}}\sum_{i=k+1}^m(e(i)-\bar{e}_n)\right|\\*
    &\qquad=o\left(d_n^{-1}\sqrt{\frac{\log n}{a_n}} \right)=o(1)\qquad P-a.s.,\\
        &\max_{k\ls m-a_n}\left|\frac{A_{k4}}{A_{k5}}\right|=\max_{k\ls m-a_n}2d_n^{-1}\frac{1}{n-m}\left|\sum_{i=1}^k(e(i)-\bar{e}_n)\right| \\*
    &\qquad=o\left(d_n^{-1}\sqrt{\frac{\log n}{n}} \right)=o(1),\qquad P-a.s.\\
\end{align*}}
Note that $A_{k5}$ is increasing in $k$ for $k\ls m$, so that
\[
\max_{k\ls m}A_{k5}=0\quad\text{and}\quad \max_{k\ls m-a_n}A_{k5}\ls -\eps \log n.
\]
Thus
\begin{align*}
    &P\left( \wm\ls m-\eps d_n^2\log n\qquad \text{inf. often} \right)\\
    &\ls P\left( \max_{k\ls m-a_n}V(k)\gs \max_{k>m-a_n}V(k)\qquad \text{inf. often} \right)\\
    &\ls P\left( \max_{k\ls m-a_n}V(k)\gs 0\qquad \text{inf. often}  \right)\\
    &\ls P\left( \max_{k\ls m-a_n}\left\{A_{k5}\left( 1+Y_n \right)\right\}\gs 0\qquad\text{inf. often} \right),\qquad\text{with }Y_n=o(1)\quad P-a.s.\\
    &=0.
\end{align*}
A similar argument gives
\[
P\left( \wm\gs m-\eps d_n^2\log n\qquad \text{inf. often} \right)=0.
\]
Hence assertion a) is proven.

For assertion b) we first need somewhat stronger bounds for the above sums, but only in a $P$-stochastic sense. Theorem B.3 in Kirch~\cite{kirchphd} gives a \haj-R\'{e}nyi type inequality if certain moment conditions of the sums are fulfilled. This yields here ($C>0$ arbitrary fixed constant)
\begin{align*}
	&\max_{1\ls k\ls n}\left|\frac{1}{\sqrt{k}}\sum_{j=1}^n (e(j)-\bar{e}_n)\right|=o_P(\sqrt{\log n}),\\
	&\max_{1\ls k\ls n}\left|\frac{1}{\sqrt{n}}\sum_{j=1}^n (e(j)-\bar{e}_n)\right|=O_P(1),\\
	&|d_n^{-1}|C^{1/2}\max_{k\ls m-C d_n^{-2}}\left|\frac{1}{m-k}\sum_{j=k+1}^m(e(j)-\bar{e}_n)\right|=O_P(1),
\end{align*}
where the last line follows because for all $\gamma_n\gs 1$ (in particular for $\gamma_n\to\infty$), $D>0$ and some $\delta>0$
\begin{align*}
	&P\left(\gamma_n^{1/2}\max_{\gamma_n\ls k\ls n}\frac 1 k \sum_{i=1}^k(e(i)-\bar{e}_n)\gs D  \right)\\
	&\ll \frac 1 {D^{2+\delta}} \gamma_n^{1+\delta/2}\left( \sum_{k=1}^{\gamma_n}\frac{1}{\gamma_n^{2+\delta}}k^{\delta/2}+\sum_{k=\gamma_n+1}^nk^{-2-\delta/2} \right)\\
	&\ll \frac{1}{D^{2+\delta}}.
\end{align*}
Analogously to above this yields for $b_n=C d_n^{-2}$, where $C>0$ is an arbitrary fixed constant,
\begin{align*}
&	\max_{k\ls m-b_n}\left|\frac{A_{k1}}{A_{k5}}\right|=o_P(1),\qquad 	\max_{k\ls m-b_n}\left|\frac{A_{k2}}{A_{k5}}\right|=o_P(1),\\
&	\max_{k\ls m-b_n}\left|\frac{A_{k3}}{A_{k5}}\right|=C^{-1/2}O_P(1),\qquad 	\max_{k\ls m-b_n}\left|\frac{A_{k4}}{A_{k5}}\right|=o_P(1).
\end{align*}
Similarly for $-b_n\ls k-m<0$
\begin{align*}
	&\max_{-b_n\ls k-m<0}\left|A_{k1}\right|=o_P(1),\qquad\max_{-b_n\ls k-m<0}\left|A_{k2}\right|=o_P(1),\\
	&\max_{-b_n\ls k-m<0}\left|A_{k4}\right|=o_P(1), \qquad A_{k5}=-(m-k)d_n^2+o_P(1),
\end{align*}
where the last rate is uniformly in $-b_n\ls k-m<0$. The proof can be finished analogously to the proof of Theorem~2 in Antoch et al.~\cite{anthusver95}, where we now use Theorem~1 of Section~1.5 in Doukhan~\cite{doukhan94}. 
\end{proof}

We will first formulate some auxiliary lemmas, which will enable us to prove the results in Section~\ref{sec_mainresults}.

\begin{lemma}\label{lem_bill_tri}
Let  $\xi_{n}(1),\ldots,\xi_n(n)$ be a triangular array of row-wise i.i.d. random variables
with $\E\xi_{n}(1)=0$ and $\E\xi_{n}^2(1)=\sigma^2+o(1)$ as
$n\to\infty$, then
\[
\left\{\frac{1}{\sigma\sqrt{n}}\sum_{j\ls nt}\xi_n(j): 0\ls t \ls 1\right\}\overset{D[0,1]}{\longrightarrow}\{W(t):0\ls t\ls 1\}.
\]
\end{lemma}
\begin{proof}
    It is analogous to that of Theorem 16.1 in Billingsley \cite{bill},
    since the central limit theorem holds
for triangular arrays and the proof of tightness also works analogously.
\end{proof}

\begin{lemma}\label{lemma_var}
    Assume that \eqref{LPmodel}-\eqref{eq_def_tau} with $0<\vth<1$ and let \eqref{assump_delta_2} be fullfilled.
    Moreover let assumption $(\mathcal{A})$ be fulfilled for some $0<\delta<(\nu-4)/2$, $\Delta=\nu-4-2\delta$. If additionally \eqref{assump_delta_est_2},
then
\begin{align*}
    &\var^*\left(\frac{1}{\sqrt{K}}\sum_{k=1}^Ke^*(k)\right)=\widetilde{\tau}_n^2+O\left(\sqrt{\frac{\log n}{L}}\right)
    =\tau^2+o(1)=O(1)\qquad P-a.s.,\\
&\text{where}\qquad\widetilde{\tau}_n^2=\frac 1 n
\sum_{l=0}^{n-1}\left[\frac{1}{\sqrt{K}}\sum_{k=1}^{K}e(l+k)\right]^2
\end{align*}
and $e(j)=e(j-n)$ for $j>n$.
\end{lemma}

\begin{rem}
    More careful considerations concerning $A_3$ below even yield an almost sure rate of $O\left( \frac{\log n}{L} \right)$.
\end{rem}
\begin{rem}
	'Estimator' $\widetilde{\tau}_n$ is closely related to the Bartlett window estimator with parameter $K-1$ if for this estimator one also uses a circularly extended series, precisely
	\[
	\widetilde{\tau}_n^2=\frac{1}{n}\sum_{l=1}^ne^2(l)+2\sum_{t=1}^{K-1}\left(1-\frac{t}{K}\right)\frac 1 n \sum_{l=1}^ne(l)e(l+t).
	\]
\end{rem}
\begin{proof}
    For $m<\wm$ (the other case can be dealt with in a similar way) \eqref{eq_decomp_etilde} yields the following decomposition
    \begin{align*}
    &   \var^*\left(\frac{1}{\sqrt{K}}\sum_{k=1}^Ke^*(k)\right)=A_1+A_2+2A_3,
    \end{align*}
where
\begin{align*}
    &   A_1=\frac{1}{n}\sum_{l=0}^{n-1}\left(\frac{1}{\sqrt{K}} \sum_{k=1}^Ke(l+k) \right)^2-K \bar{e}_n^2\\
    &A_2=\frac{1}{K^2L}\sum_{l=0}^{n-1}\left(\sum_{k=1}^K\left[d_n\left( 1_{\{n\gs l+k>m\}}-\frac{n-m}{n} \right)
    -\widehat{d}_n\left( 1_{\{n\gs l+k>\wm\}}-\frac{n-\wm}{n} \right)\right]\right)^2\\
    &A_3=\frac{1}{K^2L}\sum_{l=0}^{n-1}\sum_{k=1}^K(e(l+k)-\bar{e}_n)\\
    &\qquad\qquad\qquad \cdot\sum_{j=1}^K\left(d_n\left( 1_{\{n\gs l+k>m\}}-\frac{n-m}{n} \right)-\widehat{d}_n\left( 1_{\{n\gs l+k>\wm\}}-\frac{n-\wm}{n} \right)  \right).
\end{align*}
Theorem B.8 b) in Kirch~\cite{kirchphd} yields
\begin{equation*}
    A_1=\frac{1}{n}\sum_{l=0}^{n-1}\left(\frac{1}{\sqrt{K}} \sum_{k=1}^Ke(l+k) \right)^2+o\left( \frac{\log n}{L} \right)\qquad P-a.s.
\end{equation*}
Concerning $A_2$ we have
\begin{align*}
    A_2&\ll \frac{1}{L}\sum_{l=0}^{n-1}\left( d_n\frac{n-m}{n}-\widehat{d}_n\frac{n-\wm}{n} \right)^2\\
    &\quad +\frac{1}{K^2L}\sum_{l=0}^{n-1}\left( \sum_{k=1}^K\left( d_n 1_{\{n\gs l+k>m\}}-\widehat{d}_n1_{\{n\gs l+k>\wm\}} \right) \right)^2\\
    &=A_{21}+A_{22},
\end{align*}
now Theorem~\ref{th_asym_limit}~a) and \eqref{assump_delta_est_2} yield
\begin{align*}
    A_{21}&\ll K\left( d_n^2\frac{(m-\wm)^2}{n^2}+(d_n-\widehat{d}_n)^2 \right)\\
    &\ll K d_n^2 \frac{|m-\wm|}{n}+(d_n-\widehat{d}_n)^2 K=O\left( \frac{\log n}{L} \right)\qquad P-a.s.
\end{align*}
Similar arguments give
\begin{align*}
    A_{22}&=\frac{1}{K^2L}\sum_{l=0}^{n-1}\left( \sum_{k=1}^Kd_n 1_{\{m<l+k<\wm\}} \right)^2+\frac{1}{K^2L}\sum_{l=0}^{n-1}\left( \sum_{k=1}^K(d_n-\widehat{d}_n)1_{\{n\gs l+k>\wm\}} \right)^2\\
    &\ll \frac 1 L d_n^2 |m-\wm| + K (d_n-\widehat{d}_n)^2 =O\left( \frac{\log n}{L} \right)\qquad P-a.s.
\end{align*}

Note that
\begin{align*}
	&\E\left( \frac{1}{\sqrt{K}}\sum_{j=1}^Ke(j)
	\right)^2=\var(e(0))+2\frac{1}{K}\sum_{i<j}\cov(e(i),e(j))\\
	&=\var(e(0))+2\sum_{j=1}^{K-1}\cov(e(0),e(j))\left( 1-\frac{j}{K} \right)\\
	&=\tau^2+o(1),
	\end{align*}
	where the last line follows because of the absolute summability of the covariance function.
Thus, $A_1=\tau^2+o(1)=O(1)$ $P-a.s.$ because of (3.6.7) in Kirch~\cite{kirchphd}. We do not need assumption (3.6.6) there because (3.6.9) can be strengthened using the Minkoswki inequality as follows:
\begin{align*}
    &   P\left( \left|\frac{1}{n}\sum_{i=0}^{n-1}Y(i)\right|\gs \eps \right)\ls P\left( \left|\frac{1}{K}\sum_{k=1}^K\frac{1}{L}\sum_{l=0}^{L-1}Y(Kl+k)\right|\gs \eps \right)\\
    &\ll \E \left|\frac{1}{K}\sum_{k=1}^K\frac{1}{L}\sum_{l=0}^{L}Y(Kl+k)\right|^{2+\delta}\\
    &\ls \frac{1}{K^{2+\delta}}\left[ \sum_{k=1}^K\left( \E\left| \frac{1}{L}\sum_{l=0}^{L-1}Y(Kl+k) \right|^{2+\delta} \right)^{1/(2+\delta)} \right]^{2+\delta}\\
    &\ll  L^{-(2+\delta)/2}.
\end{align*}

The Cauchy-Schwartz inequality yields
\[
A_3\ls \sqrt{A_1\,A_2}=O\left( \sqrt{\frac{\log n}{L}} \right)\qquad P-a.s.
\]
Putting everything together we arrive at the assertion.
\end{proof}

\begin{lemma}\label{lem_diff_block}
    Assume that \eqref{LPmodel}-\eqref{eq_def_tau} with $0<\vth<1$ and let \eqref{assump_delta_2} and \eqref{assump_delta_est_2} be fullfilled.
Moreover let assumption $(\mathcal{A})$ be fulfilled for some $\delta,\Delta >0$.
Then,
{\allowdisplaybreaks\begin{align*}
&\max_{\substack{1\ls l<L\\1\ls k\ls
K}}\frac{1}{\sqrt{Kl+1}}\left|\sum_{j=k+1}^Ke^*(Kl+j)\right|\\*
&\qquad=
\max_{\substack{1\ls l<L\\1\ls k\ls
K}}\frac{1}{\sqrt{Kl+1}}\left|\sum_{j=k+1}^K(e(U(l)+j)-\bar{e}_n)\right|+R^{(1)}_n(U),\\
&\max_{1\ls k\ls K}\left|\frac{1}{\sqrt{k}}\sum_{i=1}^ke^*(i)\right|=
\max_{1\ls k\ls K}\left|\frac{1}{\sqrt{k}}\sum_{i=1}^k[e(U(0)+i)-\bar{e}_n]\right|+R^{(2)}_n(U),\\
&\max_{1\ls k<K}\left|\frac{1}{\sqrt{K-k}}\sum_{i=k+1}^Ke^*(i)\right|\\
&\qquad=
\max_{1\ls k<K}\left|\frac{1}{\sqrt{K-k}}\sum_{i=k+1}^K[e(U(0)+i)-\bar{e}_n]\right|+R^{(3)}_n(U),\\
&\max_{\substack{1\ls l<L\\1\ls k\ls
K}}\frac{1}{\sqrt{KL}}\left|\sum_{j=k+1}^Ke^*(Kl+j)\right|\\*
&\qquad=
\max_{\substack{1\ls l<L\\1\ls k\ls
K}}\frac{1}{\sqrt{KL}}\left|\sum_{j=k+1}^K(e(U(l)+j)-\bar{e}_n)\right|+R^{(4)}_n(U),
\end{align*}}
with $|R_n^{(i)}(U)|\ls R_n, i=1,2,3,$ where $R_n$ does not depend on $\{U(\cdot)\}$ and
\[
R_n=o\left( \sqrt{\log n} \right)\qquad P-a.s.,
\]
and $R_n^{(4)}(U)\ls \widetilde{R}_n$, where $\widetilde{R}_n$ does not depend on $\{U(\cdot)\}$
 and $\widetilde{R}_n=o(1)$ $P-a.s.$
\\If additionally $K d_n^2=o(1)$, then  $R_n=o(1)$ $P-a.s.$
\end{lemma}

\begin{proof}
We  only prove the first result, the others can be proven similarly. Let $m<\wm$,
 the other case is analogous.
By \eqref{eq_decomp_etilde} we get
\begin{align*}
    &\left|\max_{\substack{1\ls l<L\\1\ls k\ls
K}}\frac{1}{\sqrt{Kl+1}}\left|\sum_{j=k+1}^Ke^*(Kl+j)\right|-
\max_{\substack{1\ls l<L\\1\ls k\ls
K}}\frac{1}{\sqrt{Kl+1}}\left|\sum_{j=k+1}^K(e(U(l)+j)-\bar{e}_n)\right|\right|\\
&\ls \max_{\substack{1\ls l<L\\1\ls k\ls
K}}\frac{1}{\sqrt{Kl+1}}\Biggl|\sum_{j=k+1}^K\Bigl( d_n \left( 1_{\{n\gs U(l)+k>m\}}-\frac{n-m}{n} \right)\\
&\qquad\qquad\qquad\qquad-\widehat{d}_n\left( 1_{\{n\gs U(l)+k>\wm\}}-\frac{n-\wm}{n} \right)\Bigr)\Biggr|\\
&\ll\max_{l,k}\frac{1}{\sqrt{Kl+1}}\left|\sum_{j=k+1}^K\left(d_n\frac{n-m}{n}-\widehat{d}_n\frac{n-\wm}{n}\right)\right|\\
&\qquad +\max_{l,k}\frac{1}{\sqrt{Kl+1}}\left|\sum_{j=k+1}^Kd_n 1_{\{m<U(l)+k<\wm\}}\right|\\
&\qquad +\max_{l,k}\frac{1}{\sqrt{Kl+1}}\left|\sum_{j=k+1}^K|d_n-\widehat{d}_n|1_{\{n\gs U(l)+k>\wm\}}\right|\\
&=B_1+B_2+B_3,
\end{align*}
where by Theorem~\ref{th_asym_limit}~a) and \eqref{assump_delta_est_2} we get (note that $\min(K,|m-\wm|)\ls \sqrt{K |m-\wm|}$)
{\allowdisplaybreaks
\begin{align*}
    &B_1\ll\sqrt{K}\left(|d_n|\frac{|\wm-m|}{n}+|d_n-\widehat{d}_n|\right)=O\left( \sqrt{\frac{\log n}{L}} \right)=o(1)\qquad P-a.s.,\\
    &B_2\ll \frac{1}{\sqrt{K}}\min(K,|m-\wm|)\,|d_n|\ls \min\left(\sqrt{K}|d_n|,\sqrt{|m-\wm| |d_n^2|}\right)\\*
    &=o\left( \min\left(\sqrt{\log n},\sqrt{K}d_n \right)\right)\qquad P-a.s.,\\
    &B_3\ll \sqrt{K}|d_n-\widehat{d}_n|=O\left( \sqrt{\frac{\log n}{L}} \right)=o(1)\qquad P-a.s.
\end{align*}
}This gives the assertion.
\end{proof}

\begin{lemma}
    \label{lemma_maxsums}
    Assume that \eqref{LPmodel}-\eqref{eq_def_tau} with $0<\vth<1$ and let \eqref{assump_delta_2} and \eqref{assump_delta_est_2} be fullfilled.
    Moreover let assumption $(\mathcal{A})$ be fulfilled for some $0<\delta^{(2)}+\Delta^{(2)}<(\delta^{(1)}-2)/2<(\nu-4)/2$, $\Delta^{(1)}=\nu-2-\delta^{(1)}$.
Then,
\begin{align*}
    &a)\qquad \frac{1}{\sqrt{\log n}}\max_{1\ls j\ls n}\frac{1}{\sqrt{j}}\left|\sum_{i=1}^je^*(i)\right|=O_{P^*}(1)\qquad P-a.s.\\
    &b)\qquad \max_{1\ls j\ls n}\frac{1}{\sqrt{n}}\left|\sum_{i=1}^je^*(i)\right|=O_{P^*}(1)\qquad P-a.s.\\
    &c)\qquad \text{If moreover }d_n^2 K\to 0,\text{ then}\\
    &\phantom{c)}\qquad\widehat{d}_n^{-1}C^{1/2}\max_{1\ls j\ls \wm -C\widehat{d}_n^{-2}}\frac{1}{\wm-j}\left|\sum_{i=j+1}^{\wm}(e^*(i)-\bar{e}^*_n)\right|=O_{P^*}(1)\qquad P-a.s.
\end{align*}
\end{lemma}

\begin{proof}
    The proof of Lemma~\ref{lemma_var} shows that
    \begin{equation}\label{eq_max_var}
        \E^*\left( \frac{1}{\sqrt{K}}\sum_{k=1}^K(e(U(0)+k)-\bar{e}_n) \right)^2=O(1)\qquad P-a.s.
    \end{equation}
    Analogously we get from the proof of Theorem 3.6.2 and Remark~3.6.4 in Kirch~\cite{kirchphd} for some $\delta>0$
    \begin{equation}\label{eq_max_mom}
        \begin{split}
            &\E^*\left( \max_{k=0,\ldots,K-1}\left|\frac{1}{\sqrt{K}}\sum_{j=k+1}^K(e(U(0)+j)-\bar{e}_n)\right|^{2+\delta}\right) =O(1)\qquad P-a.s.\\
            &\E^*\left(\frac{1}{({\log n})^{1+\delta/2}} \max_{1\ls k\ls K}\left|\frac{1}{\sqrt{k}}\sum_{j=1}^k(e(U(0)+j)-\bar{e}_n) \right|^{2+\delta}\right)=o(1)\qquad P-a.s.\\
    \end{split}
    \end{equation}
    Thus we get by Lemma~\ref{lem_diff_block} using the Markov resp. \haj-R\'{e}nyi inequalities
{\allowdisplaybreaks\begin{align*}
& P^*\left(\max_{j\ls
n}\frac{1}{\sqrt{j}}\left|\sum_{i=1}^je^*(i)\right|\gs c \sqrt{\log
n}\right)\\
&\ll P^*\Biggl(\max_{1\ls l <L
}\frac{1}{\sqrt{Kl+1}}\left|\sum_{i=0}^l\sum_{k=1}^Ke^*(Ki+k)\right|+
\max_{1\ls k\ls K}\frac{1}{\sqrt{k}}\left|\sum_{i=1}^k(e(U(0)+i)-\bar{e}_n)
\right|\\*
&\qquad\qquad+\max_{\substack{1\ls l<L\\1\ls k\ls
K}}\frac{1}{\sqrt{Kl+1}}\left|\sum_{j=k+1}^K(e(U(l)+j)-\bar{e}_n)\right|\gs c \sqrt{\log n}\Biggr)\\
& \ll\frac{1}{\log n}\sum_{l=1}^L\frac 1 l  \var^*\left(\frac{1}{\sqrt{K}}\sum_{k=1}^Ke^*(k)\right)\\
&\quad+\E^*\left(\frac{1}{({\log
n})^{1+\delta/2}}\max_{1\ls k\ls K}\left|\frac{1}{\sqrt{k}}\sum_{i=1}^k(e(U(0)+i)-\bar{e}_n)\right|^{2+\delta}\right)\\*
&\quad+\frac{1}{(\log n)^{1+\delta/2}}\sum_{l=1}^L\frac{1}{ l^{1+\delta/2}} E^*\left(  \max_{0\ls k <K}\left|\frac{1}{\sqrt{K}}\sum_{j=1}^k(e(U(0)+j)-\bar{e}_n)\right|^{2+\delta}  \right)\\
&=\frac 1 K \var^*\left(\sum_{k=1}^Ke^*(k)\right) +o(1)= O(1)\qquad P-a.s.,
\end{align*}}where the last line follows from Lemma~\ref{lemma_var}.
This is assertion a), assertion b) is analogous.

Concerning assertion c) we have:
\begin{align*}
& P^*\left(\widehat{d}_n^{-1}C^{1/2}\max_{j\ls
\wm-C\widehat{d}_n^{-2}}\frac{1}{\wm-j}\left|\sum_{i=j+1}^{\wm}(e^*(i)-\bar{e}^*_n)\right|\gs c\right)\\
&\ll P^*\left(\widehat{d}_n^{-1}C^{1/2}|\bar{e}_n^*|\gs c'\right)+P^*\left(\widehat{d}_n^{-1}C^{1/2}\max_{j\ls\wm-C\widehat{d}_n^{-2}}
\frac{1}{\wm-j}\left|\sum_{l=\lceil j/K\rceil}^{\lfloor \wm/K\rfloor -1}\sum_{k=1}^Ke^*(Kl+k)\right|\gs c'\right)\\
&\qquad+P^*\left(\widehat{d}_n^{-1}C^{1/2}\max_{Kl+k\ls \wm-C\widehat{d}_n^{-2}}\frac{1}{\wm-(Kl+k)}\left|\sum_{j=k+1}^Ke^*(Kl+j)\right|\gs c'\right)\\
&\qquad +P^*\left(\widehat{d}_n^{-1}C^{1/2}\max_{j\ls\wm-C\widehat{d}_n^{-2}}\frac{1}{\wm-j}\left|\sum_{i=j+1}^{K\lceil j/K\rceil}e^*(i)\right|\gs c'\right)\\
&\qquad+P^*\left(\widehat{d}_n C^{-1/2}\left|\sum_{j=K\lfloor \wm/K\rfloor +1}^{\wm}e^*(j)\right|\gs c'\right)\\
&=D_1+D_2+D_3+D_4+D_5.
\end{align*}
The Chebyshev inequality, Lemma~\ref{lemma_var}, and the fact that $\widehat{d}_n^{-1}n^{-1/2}=o(1)$ $P-a.s.$ yield
\begin{align*}
&D_1=P^*\left(\widehat{d}_n^{-1}C^{1/2}|\bar{e}_n^*|\gs c'\right)
\ll C \widehat{d}_n^{-2}n^{-1} \frac{1}{K}\var^*\left(\sum_{i=1}^Ke^*(i)\right)=o(1)\qquad P-a.s.
\end{align*}
Moreover the \haj-R\'{e}nyi inequality gives (note that $\sum_{j=a_n}^{b_n}1/j^2=O(1/a_n)$)
\begin{align*}
D_2
& \ll C \widehat{d}_n^{-2}\var^*\left(\sum_{k=1}^Ke^*(k)\right)\left[\frac 1 {K^2}\sum_{l=C\widehat{d}_n^{-2}/K+1}^{\lfloor\wm/K\rfloor-1}\frac{1}{(l-1)^2} +\frac  {\widehat{d}^2_n}{C K} \right]\\
&\ll \frac 1 K\var^*\left(\sum_{i=1}^Ke^*(i)\right)=O(1) \qquad P-a.s.
\end{align*}
Since $K\widehat{d}^2_n\to 0$ and by Lemma~\ref{lem_diff_block} we get
\begin{align*}
D_3
&\ll \sum_{l=1}^{(\wm-C\widehat{d}_n^{-2})/K}P^*\left(\frac{C^{1/2}}{\widehat{d}_n}\frac{1}{\wm-Kl-K}\left|\max_{1\ls k\ls K}\sum_{j=k+1}^Ke^*(Kl+j)\right|\gs c''\right)\\
& \ll\left(\sqrt{K} \widehat{d}_n^{-1}C^{1/2}\right)^{2+\delta}\\
&\qquad \cdot\sum_{l=1}^{(\wm-C\widehat{d}_n^{-2})/K}\frac{1}{(\wm-Kl-K)^{2+\delta}}\E^*\left|\max_{1\ls k\ls K}\frac{1}{\sqrt{K}}\sum_{j=k+1}^K(e(U(0)+j)-\bar{e}_n)\right|^{2+\delta}\\
& \ll\left(\sqrt{K} \widehat{d}_n^{-1}C^{1/2}\right)^{2+\delta}\frac{1}{K^{2+\delta}} \sum_{j=C\widehat{d}_n^{-2}/K}^{\wm/K-2}\frac{1}{j^{2+\delta}}\\
& \ll\left(C^{-1}\widehat{d}_n^{2}K\right)^{\delta/2}=o(1)\qquad P-a.s.
\end{align*}
The arguments for $D_4$ and $D_5$ are analogous and therefore omitted.
\end{proof}

Now we are able to prove the result for the local case.

\begin{proof}[of Theorem~\ref{theorem_local}]
The proof is very close to the proof of Theorem 3 of Antoch et
al.~\cite{anthusver95}, therefore we will only sketch it.

Taking $x>0$, without loss of generality, we have for large enough $C>0$
\begin{align*}
    P^*\left( \widehat{d}_n^2(\wm^*-\wm)\ls x \right)=&P^*\left(\wm-C\widehat{d}_n^{-2} \ls \wm^*\ls \wm + x\widehat{d}_n^{-2}\right)\\
    & +P^*\left( \widehat{d}_n^2(\wm^*-\wm)\ls x, |\wm^*-\wm|>C\widehat{d}_n^{-2} \right),
\end{align*}
and the second term is smaller than
\[
P^*(\wm^*<\wm-C\widehat{d}_n^{-2})+P^*\left( \wm^*>\wm+C\widehat{d}_n^{-2} \right).
\]
Similar to the proof of Theorem~\ref{th_asym_limit} we will now show that this becomes arbitrarily small for almost all samples $X(1),\ldots,X(n)$.
First note that by Theorem~\ref{th_asym_limit}~a) there exists $\eps>0$ with
\[
0<\eps<\frac{\wm}{n}<1-\eps\qquad P-a.s.
\]
We have an analogous decomposition $V_j^*=A_{j1}^*+ A_{j2}^*+A_{j3}^*+A_{j4}^*+A_{j5}^*$ as in \eqref{eq_decomp_V}, where we replace $m$ by $\wm$ and $e(\cdot)$ by $e^*(\cdot)$.

Then, by Lemma \ref{lemma_maxsums},  
\begin{align*}
&\max_{j\ls
\wm-C\widehat{d}_n^{-2}}\left|\frac{A^*_{j1}}{A_{j5}^*}\right|
\ll \widehat{d}_n^{-2} n^{-1} \max_{j\ls
n}\frac{1}{j}\left|\sum_{i=1}^je^*(i)\right|^2
=O_{P^*}\left( \widehat{d}_n^{-2}n^{-1} \log n \right)\\
&=o_{P^*}(1)\qquad P-a.s.,
\end{align*}
since $\widehat{d}_n^{-1}n^{-1/2}(\log n)^{1/2}\to 0$ $P-a.s.$
Similarly we get
{\allowdisplaybreaks \begin{align*}
    &\max_{j\ls
\wm-C\widehat{d}_n^{-2}}\left|\frac{A^*_{j2}}{A_{j5}^*}\right|
\ll n^{-1/2} \widehat{d}_n^{-2}\max_{j\ls
\wm-C\widehat{d}_n^{-2}}\frac{1}{\wm-j}\left|\sum_{i=j+1}^{\wm}(e^*(i)-\bar{e}_n^*)\right|\max_{1\ls j\ls n}\frac{1}{\sqrt{n}}\left|\sum_{i=1}^je^*(i)\right|\\*
&\quad=O_{P^*}\left( \widehat{d}_n^{-1}n^{-1/2} \right)=o_{P^*}(1)\qquad P-a.s.\\
    &\max_{j\ls
\wm-C\widehat{d}_n^{-2}}\left|\frac{A^*_{j3}}{A_{j5}^*}\right|
\ls \widehat{d}_n^{-1}\max_{j\ls
\wm-C\widehat{d}_n^{-2}}\frac{1}{\wm-j}\left|\sum_{i=j+1}^{\wm}(e^*(i)-\bar{e}_n^*)\right|\\*
&\quad=C^{-1/2}O_{P^*}(1)\qquad P-a.s.\\
    &\max_{j\ls
    \wm-C\widehat{d}_n^{-2}}\left|\frac{A^*_{j4}}{A_{j5}^*}\right|\ll\widehat{d}_n^{-1}n^{-1/2}\max_{1\ls j\ls n}\frac{1}{\sqrt{n}}\left|\sum_{i=1}^je^*(i)\right|=O_{P^*}(\widehat{d}_n^{-1}n^{-1/2})\\*
    &\quad=o_{P^*}(1)\qquad P-a.s.
\end{align*}}This gives as in the proof of Theorem 1 in Antoch et al.~\cite{anthusver95} (and similar to the proof of Theorem~\ref{th_asym_limit}) that $P^*( \wm^*<\wm-C\widehat{d}_n^{-2} )$ becomes arbitrarily small $P-a.s.$ Similar arguments hold true for $P^*(\wm^*>\wm+C\widehat{d}_n^{-2})$.

So it suffices to consider ($C$ large enough)
\[
P^*(\wm-C\widehat{d}_n^{-2}\ls \wm^*\ls \wm+x\widehat{d}_n^{-2})=P^*\left( \max_{(j-\wm)\widehat{d}_n^2\in [-C,x]}V_j^*\gs\max_{(j-\wm)\widehat{d}_n^2\not\in [-C,x]}V_j^* \right).
\]

For $-C\widehat{d}_n^{-2}\ls j-\wm<0$ we get similarly to above
\begin{align*}
	&\max_{-C\widehat{d}_n^{-2}\ls j-\wm<0}|A_{j1}^*|=O_{P^*}\left( \widehat{d}_n^{-2}n^{-1} \right)=o_{P^*}(1)\qquad P-a.s.\\
	&\max_{-C\widehat{d}_n^{-2}\ls j-\wm<0}|A_{j2}^*|=O_{P^*}\left(\widehat{d}_n^{-1} n^{-1/2} (\log n)^{-1/2} \right)=o_{P^*}(1)\qquad P-a.s.\\
	&\max_{-C\widehat{d}_n^{-2}\ls j-\wm<0}|A_{j4}^*|=O_{P^*}\left( \widehat{d}_n^{-1}n^{-1/2} \right)=o_{P^*}(1)\qquad P-a.s.
\end{align*}
Analogous arguments give the assertion for $j>\wm$.
We can now essentially finish the proof as in Antoch et al.
\cite{anthusver95}. Note
that $P-a.s.$ we have
\[
V_j^*=
\begin{cases}
0,&j=\wm,\\[2mm]
2\widehat{d}_n\,\sum_{i=j+1}^{\wm}e^*(i)-\widehat{d}_n^2(\wm-j)+o_{P^*}(1)
,&\wm-C\widehat{d}_n^{-2}\ls j<\wm,\\[2mm]
-2\widehat{d}_n\sum_{i=\wm+1}^{j}e^*(i)+\widehat{d}_n^2(\wm-j)+o_{P^*}(1)
,&\wm<j\ls \wm-x\widehat{d}_n^{-2},
\end{cases}
\]
A similar argument as in Lemma~\ref{lemma_maxsums} yields
\[
\frac{1}{\sqrt{K}}\max_{1\ls k\ls K}\left|\sum_{i=k+1}^Ke^*(i)\right|=O_{P^*}(1)\qquad P-a.s.
\]
Thus we get for $\wm-C\widehat{d}_n^{-2}\ls j<\wm$
\[
V_j^*=2\widehat{d}_n\sqrt{K}\sum_{l=\lceil j/K\rceil}^{\lfloor \wm/K\rfloor}\frac{1}{\sqrt{K}}\sum_{k=1}^Ke^*(Kl+k)-\widehat{d}_n^2(\wm-j)+o_{P^*}(1)
,
\]
and a similar equation in the other case.
Consider the process
\[
\widetilde{V}_n^*(s)=
\begin{cases}
0,&s=0,\\[4mm]
\widehat{d}_n\sqrt{K}\sum_{l=\lceil (\wm-\widehat{d}_n^{-2}s\tau^2)/K\rceil}^{\lfloor \wm/K\rfloor}\frac{1}{\sqrt{\tau K}}\sum_{k=1}^Ke^*(Kl+k)+s/2
,&s<0,\\[4mm]
-\widehat{d}_n\sum_{l=\lceil(\wm+1)/K\rceil}^{\lfloor (\wm+\widehat{d}_n^{-2}\tau^2s)/K\rfloor}\frac{1}{\sqrt{\tau K}}\sum_{k=1}^Ke^*(Kl+k)-s/2
,&s>0.\\
\end{cases}
\]
Asymptotically the maximum of $\{V_j^*\}_j$ is the same as the maximum over $\{\widetilde{V}_n^*(s)\}_{s\in\mathbb{R}}$.  By Lemmas~\ref{lem_bill_tri} and \ref{lemma_var} the process $\{\widetilde{V}_n^*(s)\}_{s\in\mathbb{R}}$ converges $P-a.s.$ to the process $\{W(s)-|s|/2:s\in\mathbb{R}\}$.
Hence, as $L\to\infty$,
\begin{align*}
&   P^*\left( \max_{(j-\wm)\widehat{d}_n^2\in[-C,x]} V_j^*\gs \max_{(j-\wm)\widehat{d}_n^2\not\in[-C,x]} V_j^*\right)\\
&\to P^*\left( \sup_{-C\ls s\ls x} (W(s)-|s|/2)\gs \sup_{s\in ]x,C]}(W(s)-|s|/2) \right)\\
&\qquad=P^*\left(-C\ls \arg\max_s(W(s)-|s|/2)\ls x\right)\qquad P-a.s.
\end{align*}
Letting $C\to\infty$ yields the desired result.
\end{proof}

\bibliographystyle{plain}
\bibliography{./literatur}
\end{document}